\UseRawInputEncoding
\documentclass[11pt]{article}
\usepackage{amsfonts}
\usepackage{mathrsfs}
\usepackage{amsmath,amsfonts,mathrsfs,amssymb,color}
\usepackage{indentfirst}
\usepackage{graphicx}
\usepackage{float}
\usepackage{subfigure}

\numberwithin{equation}{section}

\newtheorem{theo.}{\quad\, Theorem}[section]
\newtheorem{defi.}{\quad\, Definition}[section]
\newtheorem{lemm.}{\quad\, Lemma}[section]
\newtheorem{coro.}{\quad\, Corollary}[section]

\begin{document}

\title{Existence, uniqueness, and approximability of solutions to the classical Melan equation in suspension bridges$^*$ }
\author{Jinxiang Wang$^{\ast}$
\\
 {\small  Department of Applied Mathematics, Lanzhou University of Technology, Lanzhou, P.R. Chin}\\
}
\date{} \maketitle
\footnote[0]{E-mail address: wjx19860420@163.com(Jinxiang Wang),  
} \footnote[0] {$^*$Corresponding author: Jinxiang Wang.\   }

\begin{abstract}
\baselineskip 18pt
The classical Melan equation modeling suspension bridges is considered. We first study the explicit expression and the uniform positivity of the analytical solution for the simplified ``less stiff'' model, based on which we develop a monotone iterative technique of lower and upper solutions to investigate the existence, uniqueness and approximability of the solution for the original classical Melan equation. 
The applicability and the efficiency of the monotone iterative technique for engineering design calculations are discussed by verifying some examples of actual bridges. Some open problems are suggested.

 \end{abstract}

\vskip 3mm
{\small\bf Keywords.} {\small Melan equation, suspension bridges, nonlinear nonlocal term, fourth order boundary value problem,  monotone iterative technique.}

\vskip 3mm

{\small\bf MR(2000)\ \ \ 34B10, \ 34B18}

\baselineskip 22pt

\section{Introduction and historical overview}

The theory of suspension bridges has been in existence for nearly 200 years, and a historical overview of the three most important stages of development on the suspension bridge design theory can be found in the fascinating survey of Bounopane and Billington [1]. Let us briefly outline the history here. From the very beginning, the core concern of suspension bridge theory is the vertical deflection of the bridge deck under an applied live load, and the corresponding mathematical models usually are one-dimensional in nature. 

 In 1823,  Navier put forward the theory of the ``unstiffened bridge deck'' in his celebrated report [2], which is perhaps the earliest treatise for suspension bridges model containing mathematical analysis. Navier considered that any loads applied to the deck of a bridge will be totally supported by the stiffened cable, that is, the bridge deck can be extremely flexible and does not require any stiffening effects. In [2], some second-order ordinary differential equations are derived to calculate the deflection of the structure, and suggestions for the planning of suspension bridges are also proposed as a theoretical application. However, it was gradually observed that unstiffened suspension bridges in accordance with Navier's theory are prone to exhibit large amplitude wind-induced oscillations, 
 which suggests that the model introduced by Navier is oversimplified.

With the application of truss and other vertical stiffening methods to the bridge deck, 
  the elastic theory of the stiffened suspension bridge was first introduced by Rankine [3] in 1858. Rankine considered that the displacement of the bridge under live load was influenced most by the vertical deck stiffness, and treated the stiffened span essentially as a simply supported Euler-Bernoulli beam. Although stiffened suspension bridges that conform to Rankine's elastic theory showed less wind-induced vertical undulations, to build a long span bridge according to this model could be prohibitively expensive. Moreover, related experiments showed that there exists great errors in predicting the deformation of long span stiffened suspension bridges by using elastic theory.

At the end of the 19th century, the so-called deflection theory for suspension bridges was originally developed by the Austrian engineer Josef Melan in monograph [4] whose first edition goes back to 1888. By accounting for the fact that the deflection of the bridge deck caused by live load will result in an additional tension and shape change in the cable, the deflection theory reintroduces the important effect of cable stiffness and shows that the overall structure could support larger live loads with less stiff bridge decks. The economy that could be achieved with the deflection theory makes long span bridges practical, and the first suspension bridge designed with the deflection theory was the famous Manhattan Bridge built in 1909. In 1929, Steinman included the deflection theory in the second edition of his classic Practical Treatise [5]. Since then, more and more world-famous long-span suspension bridges have been designed and built under the guidance of the deflection theory, which has gradually become a widely recognized classical theory.

Melan's deflection theory is acknowledged as a milestone contribution to the understanding of suspension bridges. Today, although the engineering design of suspension bridges is often performed using finite displacement theory (see [6-7]) based on computer, 
the deflection theory is still used as the main analytical method in the preliminary design and practical calculation of suspension bridges (see, e.g., [8-20]), which enables designers to understand the influence of key parameters and the basic behavior of structure quickly, and to validate the complex simulations.

In deflection theory, a suspension bridge is modeled as an elastic beam (the deck) suspended to a sustaining cable (see Figure 1 below).
\begin{figure}[H]
\centering
\includegraphics[width=1.0\linewidth]{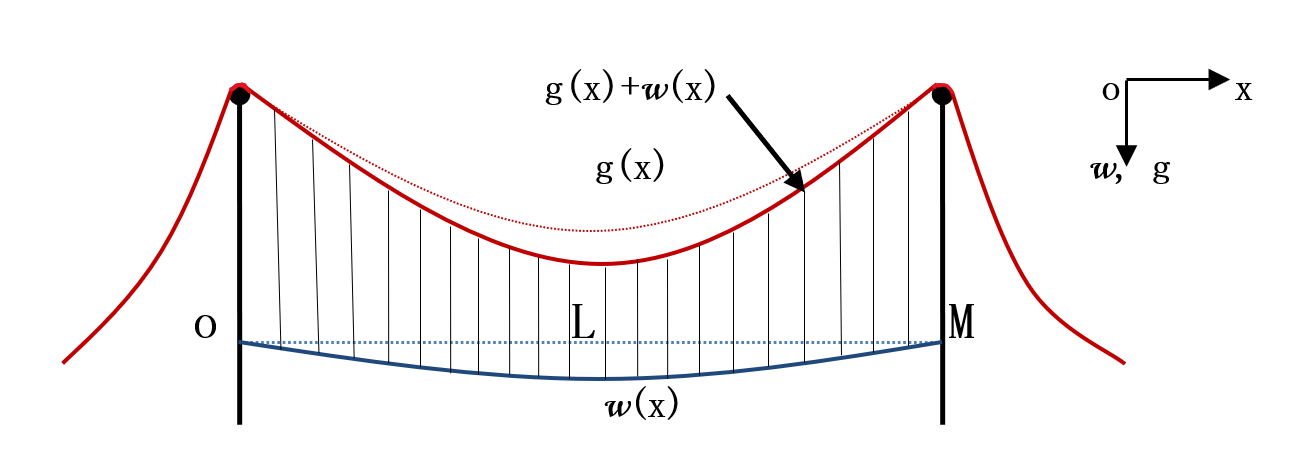}
\caption{ Beam (blue) sustained by a cable (red) through
parallel hangers.}
\end{figure}
\noindent To interpret the core problems such as the moments and shears in the suspension bridge, Melan [4, p.77] suggested a fourth-order equation to describe the vertical deflection of the bridge deck and his equation reads:
$$
EI w''''(x)-(H+h(w))w''(x)+\frac{q}{H}h(w)=p(x),\ \ \ \ x\in (0,L),
\eqno (1.1)
$$
where $w=w(x)$ denotes the vertical displacement of the beam representing the deck (positive if directed downward), $L$ is the distance between the two towers and it is also the span of the beam hinged at the two endpoints which means that the boundary conditions to be associated with (1.1) are $$
w(0)=w(L)=w''(0)=w''(L)=0.\eqno (1.2)
$$
$E$ and $I$ in (1.1) are respectively the elastic modulus of the material composing the deck and the moment of inertia of the cross section so that $EI$ is the flexural rigidity, $H$ is the horizontal tension of the cable when subject to the dead load $q$, and $h(w)$ is a nonlocal term under the integral sign (which will be detailed below in Section 8) representing the additional tension in the cable due to the live load $p=p(x)$. 

Melan's equation (1.1) is called {\it the fundamental equation of the theory of suspension bridge} by Von Karman and Biot in the monograph [21], in which a specific approximate expression of the nonlocal term $h(w)$ is also given as follows:
$$
h(w)=\frac{E_{c}A_{c}}{L_{c}}\frac{q}{H}\int_0^L w(x)dx,\eqno (1.3)
$$
where $A_{c}$ denotes the cross-sectional area of the cable, $E_{c}$ denotes the modulus of elasticity of the material, and $L_{c}$ is the length of the cable subject to the dead load $q$. The above approximate formula and the corresponding specific form of Melan equation (substituting (1.3) into (1.1)):
$$
EI w''''(x)-(H+\frac{E_{c}A_{c}}{L_{c}}\frac{q}{H}\int_0^L w(x)dx)w''(x)+\frac{q^{2}}{H^{2}}\frac{E_{c}A_{c}}{L_{c}}\int_0^L w(x)dx=p(x),\ \ \ \ x\in (0,L)\eqno (1.4)
$$
given by [21]  (see Appendix in Section 8 for detailed derivation of the equation)  
have been used by many engineering works and textbooks since then due to their conciseness. 

From a mathematical point of view, the term $h(w)$ depending on the deflections $w$ with determination (1.3) 
 makes (1.1) a nonlinear nonlocal equation, which brings great acknowledged difficulties to the qualitative research of analytical solutions. For this reason, since the deflection theory of suspension bridges was
established more than 100 years ago, several simplifications and approximate methods have been developed in engineering literature for analysing and solving the fundamental equations (1.1) and (1.4). One of the most commonly used approximate solution methods for (1.1) and (1.4) is the trial-and-error method (see Steinmann [5]): assuming a certain value for $h$, and with this value solve (1.1), (1.2). The obtained expression for $w$ is then substituted in the integral in (1.3). The result of this substitution usually will not equal $h$ since it was taken arbitrarily, and it will be necessary to repeat the calculation with a new assumed value of $h$. These trial calculations are continued far enough to obtain $h$ and $w$, with sufficient accuracy. The trial and error method has been referred to as the ``exact method.'' The procedure of this calculation was quite time-consuming before the advent of the computer, and obviously, it lacks rigorous analysis and proof. In [21] and [22], the method of trigonometric series was proposed: the procedure first treats $h$ in (1.1) as a constant and employs trigonometric series to derive a solution expressed in terms of a sine series that incorporates $h$. Subsequently, the trial-and-error method along with interpolation technique is applied to determine the value of $h$ that satisfies the required accuracy. For other simplifications and approximate methods to fundamental equations (1.1) and (1.4) in engineering literature, refer to the method of influence lines in [23] and [24], the linearized deflection theory in [25], the substitutional beam method in [11], the gravity stiffness method in [26], and the stiffness matrix method based on Laplace transformation in [27]. As can be seen from references [5],[11], and [21-27] mentioned above, the simplifications and approximate methods in existing engineering literature demonstrate a consistent lack of rigorous analysis and proof. Specifically, the relationship between approximate and true solutions is rarely addressed, and convergence proofs are entirely absent.

In [28-29], Semper considered the ``less stiff''
case of the fundamental equation (1.1) and (1.4): If the live load is relatively small compared to the dead load, $H$ will be much greater than $h(w)$, then $H+h(w)\approx H$ 
 and the equation (1.4) can be reduced to the following:

$$
EI w''''(x)-Hw''(x)+\frac{q^{2}}{H^{2}}\frac{E_{c}A_{c}}{L_{c}}\int_0^L w(x)dx=p(x),\ \ \ \ x\in (0,L).
\eqno (1.5)
$$
Semper explored a simple finite element (h-method) approach through rigorous mathematical analysis for (1.5),(1.2) and introduced an iteration scheme to find the numerical solution of (1.5), (1.2) in [28], and then, the convergence of the iteration scheme is further discussed in [29]. In the latter part of [28] and [29], numerical experiments on some examples were also carried out to demonstrate the approach and explore the difference between (1.5) and (1.4), and the results indicate that the linear integro-differential equation (1.5) seems to be a good approximation to the nonlinear nonlocal equation (1.4), even under rather large live loads.  
In [30], Gauss-Jacobi iterative technique based on a spectral method is given to get the approximate solution of (1.5), (1.2), and the iterative results similar to [28-29] are obtained. 
For other recent studies on numerical approximation methods for (1.5), (1.2), refer to the Legendre-Galerkin spectral approximation in [31], the modified Arithmetic Mean iterative method in [32], and the Optimal Homotopic Asymptotic Method in [33].

Although numerical experiments (see, e.g., [28-30]) show that the simplified ``less stiff'' model (1.5), (1.2) can well approximate the original fundamental equation (1.4), (1.2), that is, the non-linear term in equation (1.4) does not seem to make appreciable contribution. However, as the nonlinear structural characteristics of the flexible structure of suspension bridge have been discovered and discussed (see, e.g., [34-37]), the qualitative research on the original nonlinear fundamental equation (1.4) has become an unavoidable problem.
In 2015, 
Gazzola et al. [38] paid special attention to the nonlocal term $h(w)$ in the fundamental equation (1.1) and conducted a pioneering study on the well-posedness of (1.1), (1.2) and (1.4), (1.2). In addition to the common expression (1.3), two other possible forms of $h(w)$ were proposed and the differences among them were analyzed in [38]. To prove the existence and uniqueness results of solutions for (1.1), (1.2), Gazzola et al. first discussed the a priori bound for the corresponding simple constant-coefficient linear problem, 
and consequently, under some general assumptions on $h$ the second derivative term in (1.1) was treated as a whole and then the main results were proved by using fixed point theorems when $L$ and $p$ are small; see [38, Theorem 5.1 and Theorem 5.4]. It should be emphasized that the existence and uniqueness results in [38] hold only for small solutions. A counterexample given by Gazzola et al. shows that: besides the small solution, there may exist additional large solutions, that is, the problem (1.1), (1.2) and (1.4), (1.2) appear to be ill-posed. In [38], an iterative procedure based upon the solution of the simple constant-coefficient linear problem is also proposed to approximate the solution of (1.1), (1.2) and (1.4), (1.2), some numerical experiments are carried out and the results show that:  for some values of the parameters in (1.1) and (1.4), the iterative sequence seems to have a certain monotonicity (the sequence is not monotonic, but the two subsequences of odd and even iterations appear, respectively, decreasing and increasing) and admits a unique stable fixed point; however, for some other parameters such as in the range of actual bridges engineering, the sequence appears to diverge and to be quite unstable. 
At the end of [38], Gazzola et al. proposed several open problems: Under what conditions can the well-posedness of the solution to (1.1), (1.2) and (1.4), (1.2) be ensured? Can the solution be approximated by a suitable constructive sequence? 
Under which assumptions on the parameters is the iterative scheme convergent? Are there better algorithms able to manage both the stable and unstable cases? Can these algorithms detect multiple fixed points?

From [38] and the existing engineering literature (see, e.g., previously mentioned [5],[8-27] and the references therein), we notice that,
the current approximation methods for the solution of the fundamental equations (1.1), (1.2) and (1.4), (1.2) all lack rigorous convergence proof, and the existence and uniqueness of the solution also need further research. On the other hand, the exact explicit expression and properties of the solution to simplified ``less stiff'' model (1.5), (1.2) have not been addressed in numerical analysis works [28-33], and the relationship between the simplified model and the original fundamental equation (1.4), (1.2) also needs further detailed analysis. Motivated by the above two factors, 
the aim of this paper is to make further research on the well-posedness and approximation of solution for the fundamental equation (1.4), (1.2) based on the discussion of the analytical solution of simplified ``less stiff'' model (1.5), (1.2). 
Concretely, we first derive the exact explicit expression for the unique analytical solution of simplified ``less stiff'' model (1.5), (1.2) which provides a theoretical cornerstone for understanding the essence of the simplified model and serves as an exact benchmark, offering the ultimate criterion for evaluating and validating the performance of the aforementioned numerical solution work [28-33]. 
Then we prove the uniform positivity of the analytical solution within a certain range of parameters and establish a maximum principle for the corresponding integro-differential operator. 
By employing the maximum principle and analyzing the relationship between (1.5), (1.2) and (1.4), (1.2), we develop a monotone iterative technique based on the method of lower and upper solutions to 
investigate the existence, uniqueness and approximability of the solution for (1.4), (1.2). Two successively monotone iterative sequences are constructed, converging monotonically from above and below to the maximum and minimum solutions of problem (1.4), (1.2) in the sector enclosed by the lower and upper solutions, which demonstrates the algorithm's capability to detect multiple solutions of (1.4), (1.2) and provides partial answers to open problems proposed by Gazzola et al. in [38].
Besides, we prove that for parameters lying in some suitable range, the two extremal solutions are going to be equal, that is, the two monotone sequences will converge uniformly to the unique solution of (1.4), (1.2).
To the best of our knowledge, the proposed monotone iterative technique is the first approximation method for the fundamental equation (1.4), (1.2) equipped with rigorous analysis and convergence proofs. In contrast to the simplifications and approximate methods developed in the engineering literature(see, e.g., previously mentioned [5], [11], [21-27]), our constructive approximation method for the fundamental equation (1.4), (1.2) based on the analytical solution of the simplified ``less stiff'' model (1.5), (1.2) provides a way to understand the role of the nonlinear nonlocal term and the resulting nonlinear structural characteristics of the model. The applicability and efficiency of the monotone iterative technique will also be tested on some examples of actual bridges.

About the method of lower and upper solutions, as we know, for general second-order differential equation with periodic, Neumann, or Dirichlet boundary conditions, the existence of a well-ordered pair of lower and upper solutions $\alpha\leq\beta$ is sufficient to ensure the existence of a solution in the sector enclosed by them (see [39-44] and the references therein). However, it is worth noticing that, in [45] Cabada et al. pointed out that: this result no longer holds even for simple fourth-order differential equations with boundary conditions (1.2), see the counterexample in [45, Remark 3.1]. Indeed, the application of the lower and upper solutions method in boundary value problems of fourth order is heavily dependent on the conclusion of the maximum principle for the corresponding linear operators(see [45-49]). After proving the existence of a solution by using the method of lower and upper solutions, the idea of further establishing the monotone iterative technique to approximate the solution is as follows: using the maximum principle of the linear operator and starting from the ordered lower and upper solutions, one constructs two monotone sequences (one increasing and one decreasing) from the solutions of a sequence of linear problems, which converge uniformly to the extremal solutions of the original nonlinear problem between the lower and upper solutions. Therefore, the key of this technique is to select the appropriate linear operator according to the nonlinear problem and establish the corresponding maximum principle for the operator. For the specific problems we consider, after obtaining the unique analytical solution of the ``less stiff'' model (1.5), (1.2) and proving its positivity, we have actually derived the corresponding linear integro-differential operator for the original nonlinear model (1.4), (1.2) and established its maximum principle. Therefore, adopting the aforementioned monotone iterative technique becomes a natural and optimal choice. This approach simultaneously addresses the dual gaps, namely the absence of research on the analytical solution for the simplified ``less stiff'' model (1.5), (1.2) and the lack of a provably convergent approximation method for the original nonlinear model (1.4), (1.2), while also revealing the relationship between the solutions of (1.5) and (1.4).
For the maximum principle of fourth-order linear differential operator and the lower and upper solutions method of the relevant nonlinear fourth order BVPs that does not contain integral term, see
[45-58] and the references therein. As far as we know, there have been no studies on the maximum principle of fourth-order integro-differential operators and monotone iterative technique of lower and upper solutions for nonlinear nonlocal fourth order boundary value problems such as (1.4), (1.2).

The rest paper is arranged as follows: In Section 2, we study the uniqueness and explicit expression of analytical solution for ``less stiff'' model (1.5), (1.2). In Section 3, we analyze the uniform positivity of the analytical solution and establish a maximum principle for the corresponding operator. In Section 4, based upon the maximum principle 
we develop a monotone iterative technique in the presence of lower and upper solutions to investigate the existence and approximation of solutions for (1.4), (1.2), the uniqueness of the solutions for a certain range of parameters is also discussed. In Section 5, we demonstrate the application of the the main results by taking a specific form of upper and lower solutions. In Section 6, the applicability and efficiency of the iterative method for practical engineering design are discussed by verifying some examples of actual bridges. Section 7 contains our conclusions and some open problems. Finally,
for maintaining the integrity of the paper and the convenience of readers, the full derivation of the classical form of Melan equation (1.4), (1.2) will be given in Appendix in Section 8.

\section{The analytical solution for ``less stiff'' model (1.5), (1.2)}
In this Section, we first prove a uniqueness result of solutions for ``less stiff'' model (1.5), (1.2) and derive the explicit expression of the analytical solution. For simplicity, we drop some constants and consider the the following equivalent linear fourth order integro-differential equations with Navier boundary condition:
$$
\aligned
&y^{(4)}(x)-My''(x)+N\int_0^Ly(x)dx=p(x),\ \ \ \ x\in (0,L),\\
&y(0) = y(L) = y''(0) = y''(L) = 0,\\
\endaligned
\eqno (2.1)
$$
where $M=\frac{H}{EI},N=\frac{q^{2}}{EIH^{2}}\frac{E_{c}A_{c}}{L_{c}}$ are  constants, $p\in L^{1}[0,L]$.

Denote $D(K):=\{w\in W^{4,1}[0,L]:w(0)=w(L)=w''(0)=w''(L)=0\},$ 
we
define a linear operator $K:D(K)\to L^{1}[0,L]$ by
$$
Ky:=y''''-My''
$$ and investigate the Green function of $Ky=0$. To do that, let us introduce two second order linear differential operators which are associated with $K$:
$$
K_{M}y:=-y''+My,~~~D(K_{M}):=\{y\in W^{2,1}[0,L]:y(0)=y(L)=0\};
$$
$$
K_{0}y:=-y'',~~~D(K_{0}):=\{y\in W^{2,1}[0,L]:y(0)=y(L)=0\}.
$$
It is easy to verify that the Green functions of $K_{M}y=0$ and $K_{0}y=0$ are
$$G_{M}(t,s)=\left\{
\begin{aligned}
\frac{\sinh (\sqrt M t) \sinh (\sqrt M(L-s))}{\sqrt M\sinh (\sqrt M L)},\ 0\leq t\leq s\leq L,\\
\frac{\sinh (\sqrt Ms) \sinh (\sqrt M(L-t))}{\sqrt M\sinh (\sqrt M L)},\ 0\leq s\leq t\leq L
\end{aligned}
\right.\eqno (2.2)$$
and
$$G_{0}(t,s)=\left\{ \aligned
\frac{s(L-t)}{L}, \ \ \ \ \ \ \ \  \ 0\leq s\leq t\leq L,\\
\frac{t(L-s)}{L}, \ \ \ \ \ \ \ \  \ 0\leq t\leq s\leq L
\endaligned
\right.\eqno (2.3)
$$ respectively, and they have the following properties for $(t,s)\in [0,L]\times[0,L]$:   $$0\leq G_{0}(t,s)\leq G_{0}(t,t)\leq G_{0}(\frac{L}{2},\frac{L}{2})=\frac{L}{4}, \eqno (2.4)$$$$0\leq\frac{\sqrt{M}}{\sinh (\sqrt{M}L)} G_{M}(t,t)G_{M}(s,s) \leq G_{M}(t,s)\leq G_{M}(t,t)\leq G_{M}(\frac{L}{2},\frac{L}{2})=\frac{\tanh(\frac{\sqrt{M}L}{2})}{2\sqrt{M}}.\eqno (2.5)$$

Obviously,
$$Ky=K_{0}(K_{M})y,$$
and then the Green function of $Ky=0$ is
$$G(x,s):=\int^{L}_{0}G_{0}(x,t)G_{M}(t,s)dt, \ \ \ \  \ (x,s)\in [0,L]\times [0,L].\eqno (2.6)
$$
 By explicit calculation, the Green function $G(x,s)$ in (2.6) can be given by the following:
 $$G(x,s)=\left\{
\begin{aligned}
\frac{x(L-s)}{ML}-\frac{1}{M}\frac{\sinh (\sqrt M x) \sinh (\sqrt M(L-s))}{\sqrt M\sinh (\sqrt M L)},\ 0\leq x\leq s\leq L,\\
\frac{s(L-x)}{ML}-\frac{1}{M}\frac{\sinh (\sqrt M s) \sinh (\sqrt M(L-x))}{\sqrt M\sinh (\sqrt M L)},\ 0\leq s\leq x\leq L.
\end{aligned}
\right.\eqno (2.7)$$
Since $G_{M}(t,s), G_{0}(t,s)$ are nonnegative, then the Green function $G(x,s)$ given by (2.6) is nonnegative on $[0,L]\times [0,L]$.

We are now in a position to solve the problem (2.1), and the result is shown below:

\noindent{\bf Theorem 2.1}\ For any fixed constant $M>0$, assume that $$ |N|\leq\frac{1}{L( \frac{L^2}{8M})}<\frac{1}{L( \frac{L^2}{8M}+\frac{2\sinh \frac{\sqrt ML}{2}}{M^2\sinh (\sqrt ML)}-\frac{1}{M^2})}=\frac{1}{L\underset{x\in [0,L]}\max\int^{L}_{0}G(x,s)ds},\eqno (2.8)$$ then for any $p\in L^{1}[0,L]$, the nonlocal boundary value problem (2.1) has a unique solution
$$
\aligned
y(x)=&\int^{L}_{0}G(x,s)p(s)ds-
\frac{N\int^{L}_{0}G(x,s)ds}{1+N\int^{L}_{0}\int^{L}_{0}G(x,s)dsdx}\int^{L}_{0}\int^{L}_{0}G(x,s)p(s)dsdx, \ \   \ x\in [0,L].\\
\endaligned\eqno (2.9)
$$

\noindent{\bf Proof.} We first prove the existence and uniqueness of
solutions for (2.1). Observe that $y$ is a solution of (2.1) if and only if $y$ is a fixed point of the the operator $F: C[0,L]\rightarrow C[0,L]$ given by
$$
[Fy](x)=\int^{L}_{0}G(x,s)p(s)ds -N\int^{L}_{0}G(x,s)ds\int^{L}_{0}y(x)dx,
$$
where $G(x,s)$ is as in (2.7).

For $u,v\in C[0,L]$, by (2.8) we have
$$
\aligned
\|Fu-Fv\|_{\infty}&=\|N\int^{L}_{0}G(x,s)ds\int^{L}_{0}[u(t)-v(t)]dt\|_{\infty}\ \\
&\leq L\|u-v\|_{\infty}|N| \underset{x\in [0,L]}\max\int^{L}_{0}G(x,s)ds\\
&< \|u-v\|_{\infty}.\\
 \endaligned
  $$
Based on Banach fixed point theorem, there exists a  unique fixed point for the operator $F$, which assures the existence and uniqueness of
solution for (2.1).

In the sequel, we proof that the explicit expression of the unique solution of (2.1) is (2.9).

Using Picard's iterative method we know that for any $y_{0}\in C[0,L]$, the sequence given by $y_{n}=Fy_{n-1}(n\geq1)$ converges to the unique solution of (2.1). Taking
$$
y_{0}=\int^{L}_{0}G(x,s)p(s)ds,
$$
by recursive calculation we get that
$$
y_{n}(x)=y_{0}(x)+\int^{L}_{0}G(x,s)ds\int^{L}_{0}\int^{L}_{0}G(x,s)p(s)dsdxN\sum_{i=1}^{n}[(-1)^iN^{i-1}(\int^{L}_{0}\int^{L}_{0}G(x,s)dsdx)^{i-1}].
\eqno (2.10)$$
Since (2.8) can lead to $$N\int^{L}_{0}\int^{L}_{0}G(x,s)dsdx\leq|N| L \underset{x\in [0,L]}\max\int^{L}_{0}G(x,s)ds=d<1,$$
then by passing  to the limit for the Picard's iterative (2.10), we conclude that the unique solution of (2.1) is given by
$$
\aligned
y(x)&=\int^{L}_{0}G(x,s)p(s)ds-
\frac{N\int^{L}_{0}G(x,s)ds}{1+N\int^{L}_{0}\int^{L}_{0}G(x,s)dsdx}\int^{L}_{0}\int^{L}_{0}G(x,s)p(s)dsdx
, \ \   \ x\in [0,L].\\
\endaligned
 $$
 \hfill{$\Box$}

\noindent{\bf Remark 2.1} In Lemma 2.1 in [28], Semper has proved the existence and uniqueness of the solution to the ``less stiff'' model (1.5), (1.2) by using Lax-Milgram theorem. The result of [28] combining with our theorem 2.1 above imply that the condition (2.8) can be removed, that is,  the explicit expression of the analytical solution of ``less stiff'' model (1.5), (1.2) given by formula (2.9) in Theorem 2.1 is valid for any $M, N>0$. That is, we get the following conclusion:

\noindent{\bf Corollary 2.1} For any constants $M, N>0$ and $p\in L^{1}[0,L]$, the nonlocal boundary value problem (2.1) has a unique analytical solution given by (2.9).

Obviously, the 
unique analytical solution given by formula (2.9) can be seen as an improvement on numerical solution work in [28-33]. We verify and demonstrate Corollary 2.1 and (2.9) by two simple examples as follows:

\noindent{\bf Example 2.1. (see Part 6 in [28] and Part 5 in [30])} Consider problem (2.1) with $M=N=L=1$ and non-positive $p(x)=2.75-6(x+\frac{e-1}{e+1})$, that is, the following problem
$$
\aligned
&y^{(4)}(x)-y''(x)+\int_0^1y(x)dx=2.75-6(x+\frac{e-1}{e+1}),\ \ \ \ x\in (0,1),\\
&y(0) = y(1) = y''(0) = y''(1) = 0.\\
\endaligned
\eqno (2.11)
$$
Then according to Corollary 2.1 and (2.9) and by explicit calculation, problem (2.11) has the following unique non-positive solution:
$$
y(x)=5x-6\frac{\sinh x}{\sinh 1}+x^{3}
$$
with maximum deflection $\underset{x\in (0,1)}\max |y(x)|=|y(0.5)|=0.03545665.$ This result confirms that the numerical method results and the preconceived solution form in [28] (see Figures 2-3 and (6.1)) and [30] (see Table1 and Figure 2) are effective and correct.

\noindent{\bf Example 2.2. (see Part 3 in [29] and Part 5 in [30])} Consider problem (2.1) with $M=N=1$, $L=2$ and non-positive $p(x)\equiv-5$, that is, the following problem
$$
\aligned
&y^{(4)}(x)-y''(x)+\int_0^2y(x)dx=-5,\ \ \ \ x\in (0,2),\\
&y(0) = y(2) = y''(0) = y''(2) = 0.\\
\endaligned
\eqno (2.12)
$$
Then according to Corollary 2.1 and (2.9) and by explicit calculation, problem (2.12) has the following unique non-positive solution:
$$
y(x)=-\frac{15\sinh2}{6(\cosh2-1)-\sinh2}[\frac{2x-x^2}{2}+\frac{\sinh x+\sinh(2-x)-\sinh2}{\sinh2}]
$$
with maximum deflection $\underset{x\in (0,2)}\max |y(x)|=|y(1)|=0.6222.$ This result confirms that the numerical method results in [29] (see Prat 3 and Figure 1) and [30] (see Figure3 in Part 5) are effective.

\section{The positivity of the analytical solution for ``less stiff'' model (1.5), (1.2) and a Maximum principle}

In this section, we study the uniform positivity of the analytical solution (2.9) for problem (2.1), that is, we will try to establish a Maximum principle for the corresponding operator
$$
Gw=w^{(4)}-Mw''+N\int_0^Lw(t)dt  \eqno (3.1)
$$ on $D(G):=\{w\in W^{4,1}[0,L]:w(0)=w(L)=w''(0)=w''(L)=0\}.$

The result is as follows:

\noindent{\bf Theorem 3.1} For any fixed constant $M>0$, assume that $$ 0\leq N\leq \frac{4M}{L^3}\sigma(\sqrt{M}L),\eqno (3.2) $$\ where
\begin{small}
$$\sigma(\sqrt{M}L):=\frac{3(\sqrt{M}L)^{3}}{6\sqrt{M}L\cosh(\sqrt{M}L)+6\sqrt{M}L-12\sinh(\sqrt{M}L)
 -(\sqrt{M}L)^{3}},\eqno (3.3)$$
\end{small} 
 then for any $p\geq0,$ the solution $y$ of problem (2.1) given by (2.9) satisfies $y\geq0$ and $y''\leq0$; Inversely, for $p\leq0,$ the solution $y$  satisfies $y\leq0$ and $y''\geq0$.

\noindent{\bf Proof.}  Since $$
y(x)=-\int^{L}_{0}G_{0}(x,s)y''(s)ds,   \ \ \ \ \ \forall x\in(0,L),
$$ and the Green function $G_{0}(x,s)\geq0$, so to complete the proof we only need to prove the non-negative and non-positive properties of $y''$ corresponding to the two cases.

Based upon (2.6) and (2.9), by explicit calculation we have
$$
y''(x)=-\int^{L}_{0}G_{M}(x,s)p(s)ds+
\frac{N\int^{L}_{0}\int^{L}_{0}G(x,s)p(s)dsdx}{1+N\int^{L}_{0}\int^{L}_{0}G(x,s)dsdx}\int^{L}_{0}G_{M}(x,s)ds.
\eqno (3.4)$$
Now, we first prove that if (3.2) hold, then $y''$ given above is nonpositive for any $p\geq0$, that is
$$
-\int^{L}_{0}G_{M}(x,s)p(s)ds+
\frac{N\int^{L}_{0}\int^{L}_{0}G(x,s)p(s)dsdx}{1+N\int^{L}_{0}\int^{L}_{0}G(x,s)dsdx}\int^{L}_{0}G_{M}(x,s)ds\leq 0.
\eqno (3.5)$$
Since the Green function $G(x,s)$ is nonnegative, (3.5) is equivalent to the following inequality:
$$
\aligned
&
N\int^{L}_{0}\int^{L}_{0}G(x,s)p(s)dsdx\int^{L}_{0}G_{M}(x,s)ds\\
&\ \ \ \ \ \ \ \ \ \ \ \ \ \ \ \ \ \  \ \ \ \ \ \ \ \ \ \ \ \  \leq (N\int^{L}_{0}\int^{L}_{0}G(x,s)dsdx+1)\int^{L}_{0}G_{M}(x,s)p(s)ds.\\
\endaligned
\eqno (3.6)
$$
We will then explicitly calculate both sides of the inequality (3.6) in detail.

By (2.7), we have
$$
 \aligned
&\int^{L}_{0}\int^{L}_{0}G(x,s)dsdx
=\int^{L}_{0}(\frac{Lx-x^2}{2M}+\frac{\sinh(\sqrt{M}x)+\sinh(\sqrt{M}(L-x))-\sinh(\sqrt{M}L)}{M^2\sinh(\sqrt{M}L)})dx \\
&\ \ \ \ \ \ \ \ \ \ \ \ \ \ \ \ \ \ \ \ \ \ \ \ \ \ \ =\frac{(\sqrt{M}L)^3\sinh(\sqrt{M}L)+24(\cosh(\sqrt{M}L)-1)-12\sqrt{M}L\sinh(\sqrt{M}L)}{12M^2\sqrt{M}\sinh(\sqrt{M}L)}\\
&\ \ \ \ \ \ \ \ \ \ \ \ \ \ \ \ \ \ \ \ \ \ \ \ \ \ \ :=\zeta_{\sqrt{M},L},\\ \endaligned\eqno (3.7)
 $$
this combine with (2.2) can deduce that the right-hand side of (3.6) is as following:
\begin{small}
$$
 \aligned
&(N\int^{L}_{0}\int^{L}_{0}G(x,s)dsdx+1)\int^{L}_{0}G_{M}(x,s)p(s)ds= \\
&\frac{(N\zeta+1)\sinh (\sqrt M(L-x))}{\sqrt M\sinh (\sqrt M L)}\int^{x}_{0}\sinh (\sqrt Ms)p(s)ds+\frac{(N\zeta+1)\sinh (\sqrt Mx)}{\sqrt M\sinh (\sqrt M L)}\int^{L}_{x}\sinh (\sqrt M(L-s))p(s)ds;\\
\endaligned\eqno (3.8)
$$
\end{small}
Based on (2.7) and exchange the order of integration, we have
\begin{small}
$$
 \aligned
&\int^{L}_{0}\int^{L}_{0}G(x,s)p(s)dsdx \\
&=\int^{L}_{0}\{[\frac{L-x}{ML}\int^{x}_{0}sp(s)ds-\frac{\sinh (\sqrt M(L-x))}{M\sqrt M\sinh (\sqrt M L)}\int^{x}_{0}\sinh (\sqrt M s)p(s)ds]\\
&\ \ \ \ \ \ \ \ \ \ \ \ \ \ \ \ \ \ \ +[\frac{x}{ML}\int^{L}_{x}(L-s)p(s)ds-\frac{\sinh (\sqrt Mx)}{M\sqrt M\sinh (\sqrt M L)}\int^{L}_{x}\sinh (\sqrt M (L-s))p(s)ds]\}dx\\
&=\frac{1}{2M}\int^{L}_{0}s(L-s)p(s)ds+\frac{1}{M^2\sinh (\sqrt M L)}\int^{L}_{0}[\sinh (\sqrt M s)+\sinh (\sqrt M(L-s))-\sinh (\sqrt M L)]p(s)ds \\
&=\frac{1}{2M}\int^{L}_{0}s(L-s)p(s)ds-\frac{1}{M^2\sinh (\sqrt M L)}\int^{L}_{0}\theta(s)p(s)ds, \\
\endaligned\eqno (3.9)
$$
\end{small}
where $$\theta(s)=\sinh (\sqrt ML)-\sinh (\sqrt M(L-s))-\sinh (\sqrt Ms).\eqno (3.10)$$
By (2.2) we have
$$
 \aligned
&\int^{L}_{0}G_{M}(x,s)ds \\
&=\frac{\sinh (\sqrt M(L-x))}{\sqrt M\sinh (\sqrt M L)}\int^{x}_{0}\sinh (\sqrt Ms)ds+\frac{\sinh (\sqrt Mx)}{\sqrt M\sinh (\sqrt M L)}\int^{L}_{x}\sinh (\sqrt M(L-s))ds\\
&=\frac{\sinh (\sqrt ML)-\sinh (\sqrt M(L-x))-\sinh (\sqrt Mx)}{ M\sinh (\sqrt M L)}=\frac{1}{ M\sinh (\sqrt M L)}\theta(x),\\
\endaligned\eqno (3.11)
$$
where $\theta(x)$ is as defined in (3.10) above.

Combining (3.11) with (3.9) we can deduce that the left-hand side of (3.6) is as following:
\begin{small}
$$
 \aligned
&N\int^{L}_{0}\int^{L}_{0}G(x,s)p(s)dsdx\int^{L}_{0}G_{M}(x,s)ds \\
&=N\frac{\theta(x)}{2M^2\sinh (\sqrt M L)}\int^{L}_{0}s(L-s)p(s)ds-N\frac{\theta(x)}{M^3\sinh^2(\sqrt M L)}\int^{L}_{0}\theta(s)p(s)ds \\
&=\frac{N\theta(x)}{M^2\sinh (\sqrt M L)}\int^{x}_{0}[\frac{1}{2}s(L-s)-\frac{\theta(s)}{M\sinh(\sqrt M L)}]p(s)ds\\
&\ \ \ \ \ \ \ \ \ \ \ \ \ \ \ \ \ \ \ \ \ \ \ \ +\frac{N\theta(x)}{M^2\sinh (\sqrt M L)}\int^{L}_{x}[\frac{1}{2}s(L-s)-\frac{\theta(s)}{M\sinh(\sqrt M L)}]p(s)ds.\\
\endaligned\eqno (3.12)
$$
\end{small}
Now, according to (3.8) and (3.12), to prove the inequality (3.5) and (3.6) under condition (3.2), we only need to verify that the following two inequalities are true:
\begin{footnotesize}
$$
\frac{(N\zeta+1)\sinh (\sqrt M(L-x))}{\sqrt M\sinh (\sqrt M L)}\int^{x}_{0}\sinh (\sqrt Ms)p(s)ds\geq \frac{N\theta(x)}{M^2\sinh (\sqrt M L)}\int^{x}_{0}[\frac{1}{2}s(L-s)-\frac{\theta(s)}{M\sinh(\sqrt M L)}]p(s)ds,
\eqno (3.13)$$
\end{footnotesize}
\begin{footnotesize}
$$\frac{(N\zeta+1)\sinh (\sqrt Mx)}{\sqrt M\sinh (\sqrt M L)}\int^{L}_{x}\sinh (\sqrt M(L-s))p(s)ds\geq \frac{N\theta(x)}{M^2\sinh (\sqrt M L)}\int^{L}_{x}[\frac{1}{2}s(L-s)-\frac{\theta(s)}{M\sinh(\sqrt M L)}]p(s)ds.\eqno (3.14)$$
\end{footnotesize}
Next, we first show that condition (3.2) can guarantee (3.13) holds for $x\in[0,L].$ Denotes
$$
\xi_{\sqrt{M},L}=\frac{\sqrt M L}{2}\sinh(\sqrt M L)-\cosh(\sqrt M L)+1,
\eqno (3.15)$$
obviously, (3.13) is equivalent to
\begin{footnotesize}
$$
\frac{(N\zeta+1)\sinh (\sqrt M(L-x))}{\xi\sqrt M}\int^{x}_{0}\xi\sinh (\sqrt Ms)p(s)ds\geq \frac{N\theta(x)}{M^3\sinh (\sqrt M L)}\int^{x}_{0}[\frac{1}{2}s(L-s)M\sinh(\sqrt M L)-\theta(s)]p(s)ds.
\eqno (3.16)$$
\end{footnotesize}
Since $\xi_{\sqrt{M},L}=\underset{s\in (0,L)}\sup \{\frac{\frac{1}{2}s(L-s)M\sinh(\sqrt M L)-\theta(s)}{\sinh (\sqrt Ms)}\}$, it is easy to see that the functions in the integral terms on two sides of (3.16) satisfy
$$
\xi\sinh (\sqrt Ms)\geq \frac{1}{2}s(L-s)M\sinh(\sqrt M L)-\theta(s),  \ \ \ \ \ \forall s\in(0,L).\eqno (3.17)
$$
So, to obtain (3.16) and (3.13), we only need to verify that the following inequality holds under condition (3.2):
$$
\frac{(N\zeta+1)\sinh (\sqrt M(L-x))}{\xi\sqrt M}\geq \frac{N\theta(x)}{M^3\sinh (\sqrt M L)},  \ \ \ \ \ \forall x\in(0,L).\eqno (3.18)
$$
In fact, (3.18) is equivalent to
$$
 \aligned
\sinh (\sqrt M(L-x))&\geq N[\frac{\theta(x)\xi}{M^2\sqrt M\sinh (\sqrt M L)}-\zeta\sinh (\sqrt M(L-x)] \\
&=N\frac{\varphi(x)}{12M^2\sqrt{M}\sinh(\sqrt{M}L)},\\
\endaligned\eqno (3.19)
$$
where
\begin{small}
$$
\varphi(x)=\{12\xi[\sinh (\sqrt ML)+\sinh (\sqrt M(L-x))-\sinh (\sqrt Mx)]-(\sqrt ML)^{3}\sinh (\sqrt ML)\sinh (\sqrt M(L-x))\}.
$$
\end{small}
For any fixed $L>0,$ by analyzing the basic properties of the function $\varphi(x)$ such as monotonicity, there exists a zero point $\eta_{L}\in (0,L)$ such that $\varphi(x)<0$ for $x\in[0,\eta_{L})$ and $\varphi(x)>0$ for $x\in(\eta_{L},L)$. Since if $\varphi(x)\leq0$, then (3.19) hold for any $N>0$. So, in order for the inequality (3.19) to be true for any $x\in[0,L]$, it is only necessary for $N$ to satisfy the following condition:
$$
N\leq \underset{x\in (\eta_{L},L)}\inf \{\frac{\sinh (\sqrt M(L-x))}{\varphi(x)}\} 12M^2\sqrt{M}\sinh(\sqrt{M}L). \eqno (3.20)
$$
Since the function $\frac{\sinh (\sqrt M(L-x))}{\varphi(x)}$ is monotonically decreasing in $(\eta_{L},L)$ and $\varphi(\sqrt ML)=0$, then by L'Hopital's Rule we have
$$
 \aligned
&\underset{x\in (\eta_{L},L)}\inf \{\frac{\sinh (\sqrt M(L-x))}{\varphi(x)}\}=\lim\limits_{x\rightarrow L}\frac{\sinh (\sqrt M(L-x))}{\varphi(x)} \\
&\ \ \ \ \ \ \ =\frac{1}{6[2-2\cosh(\sqrt{M}L)+\sqrt{M}L\sinh(\sqrt{M}L)][1+\cosh(\sqrt{M}L)]
 -(\sqrt{M}L)^{3}\sinh(\sqrt{M}L)}.\\
\endaligned\eqno (3.21)
$$
Substituting (3.21) into the (3.20), we obtain the condition that $N$ needs to meet as follows:
$$
 \aligned
&N\leq \frac{ 12M^2\sqrt{M}\sinh(\sqrt{M}L)}{6[2-2\cosh(\sqrt{M}L)+\sqrt{M}L\sinh(\sqrt{M}L)][1+\cosh(\sqrt{M}L)]
 -(\sqrt{M}L)^{3}\sinh(\sqrt{M}L)} \\
&=\frac{4M}{L^3}\frac{3(\sqrt{M}L)^{3}\sinh(\sqrt{M}L)}{6[2-2\cosh(\sqrt{M}L)+\sqrt{M}L\sinh(\sqrt{M}L)][1+\cosh(\sqrt{M}L)]
 -(\sqrt{M}L)^{3}\sinh(\sqrt{M}L)}\\
&=\frac{4M}{L^3}\frac{3(\sqrt{M}L)^{3}}{6\sqrt{M}L\cosh(\sqrt{M}L)+6\sqrt{M}L-12\sinh(\sqrt{M}L)
 -(\sqrt{M}L)^{3}}\\
&=\frac{4M}{L^3}\sigma(\sqrt{M}L). \\
\endaligned\eqno (3.22)
$$
So far, we have proved that condition (3.2) can guarantee the inequality (3.13).

According to the symmetry between $\sinh (\sqrt Mx)$ and $\sinh (\sqrt M(L-x))$ in (3.13) and (3.14) and it is easy to verify that $\xi_{\sqrt{M},L}$ defined by (3.15) also satisfy $$\xi_{\sqrt{M},L}=\underset{s\in (0,L)}\sup \{\frac{\frac{1}{2}s(L-s)M\sinh(\sqrt M L)-\theta(s)}{\sinh (\sqrt M(L-s))}\},$$ then we can prove the inequality (3.14) under the condition (3.2) using similar discussion methods as (3.16)-(3.22), only by replacing two pairs of functions $\sinh (\sqrt Mx),\sinh (\sqrt M(L-x))$ and $\varphi(x),\varphi(L-x)$ respectively.
Since (3.13) and (3.14) both hold under condition (3.2), we immediately deduce that (3.5) holds, that is, we complete the proof for the case $p\geq0$.

Analogously we can prove the case $p\leq0$. \hfill{$\Box$}

\noindent{\bf Remark 3.1} Since $\frac{4M}{L^3}\sigma(\sqrt{M}L)=\frac{1}{L^5}4(\sqrt{M}L)^{2}\sigma(\sqrt{M}L)$, then the condition (3.2) in Theorem 3.1 is equivalent to$$ 0\leq N\leq \frac{4}{L^5}(\sqrt{M}L)^{2}\sigma(\sqrt{M}L).\eqno (3.23) $$

\noindent{\bf Remark 3.2} If (3.2) does not hold, although the uniform positivity for arbitrary nonnegative live loads $p$ may not be true, the solution $y_{p}$ of problem (2.1) under some specific live load $p$ may still have positivity as in Theorem 3.1.


\section{Monotone iterative technique for classical Melan equation (1.4),(1.2)}

In this section, we will develop the method of lower and upper solutions and monotone iterative technique for (1.4), (1.2). For simplicity, by assuming the flexural rigidity $EI\equiv 1$ and dropping some constants, we consider the following equivalent nonlinear nonlocal fourth order equations with Navier boundary condition (1.2):
$$
w''''(x)-(a+b\int_0^L w(x)dx)w''(x)+c\int_0^L w(x)dx=p(x),\ \ \ \ x\in (0,L),\eqno (4.1)
$$
where $a=H,b=\frac{E_{c}A_{c}}{L_{c}}\frac{q}{H},c=\frac{q^{2}}{H^{2}}\frac{E_{c}A_{c}}{L_{c}}$ are positive parameters, $p\in L^{1}[0,L]$.

We will use the following definition of lower and upper solutions.

\noindent{\bf Definition 4.1.} The function $\alpha\in W^{4,1} [0, L]$ is said to be a {\it lower solution} for the BVP (4.1), (1.2) if

$$\alpha''''(x)-(a+b\int_0^L \alpha(x)dx)\alpha''(x)+c\int_0^L \alpha(x)dx\leq p(x),\ \ \ \ x\in (0,L),
\eqno (4.2)
$$
$$
\alpha(0)= 0,\ \alpha(L) = 0,\ \alpha''(0) =0, \ \alpha''(L) = 0.
\eqno (4.3)
$$
An {\it upper solution} $\beta \in W^{4,1} [0, L]$ is defined analogously by reversing the inequalities in (4.2).

For any $ y\in W^{4,1} [0, L]$, denote \begin{small}
$$
f(x,\int_0^L y(x)dx, y''(x))=p(x)+ay''(x)+by''(x)\int_0^L y(x)dx-c\int_0^L y(x)dx,\ \ \ \ x\in (0,L).\eqno (4.4)
$$
\end{small}
The main result is as follows:

\noindent{\bf Theorem 4.1.} Assume that $M,N>0$ are constants and $N\leq \frac{4M}{L^3}\sigma(\sqrt{M}L)=\frac{4}{L^5}(\sqrt{M}L)^{2}\sigma(\sqrt{M}L)$ where $\sigma(\sqrt{M}L)$ is defined as in (3.3), the operator $G: D(G)\rightarrow L^{1}[0,L]$ is as in (3.1). If problem (4.1), (1.2) have a lower solution $\alpha$ and an upper solution $\beta$ given as in Definition 4.1 which satisfy
$$0\leq\alpha(x) \leq \beta(x), \ \ \   \ \ \ \ 0\geq\alpha''(x) \geq \beta''(x)  \ \ \ \ \ \ \ \ \text{for}\ x \in [0, L],
\eqno (4.5)
$$
and
$$ a+b\int_0^L \beta(x)dx\leq M,\ \ \   \ \ \ \ c-b\underset{x\in [0,L]}\min\beta''(x) \leq N.
\eqno (4.6)
$$
Then there exists at least one solution for (4.1), (1.2) in $[\alpha,\beta]=\{u\in C^2[0,1]| \alpha\leq u\leq\beta,\ \ \alpha''\geq u''\geq\beta''\}$. Moreover, with
$f$ defined in (4.4), the iterative sequences $\{\alpha_{n}\}$ and $\{\beta_{n}\}$ produced by the iterative procedure
$$
\aligned
&Gy_{n}(x)=f(x,\int_0^L y_{n-1}(x)dx, y_{n-1}''(x))-My_{n-1}''(x)+N\int_0^Ly_{n-1}(t)dt,\\
&y_{n}(0) = y_{n}(L) = y_{n}''(0) = y_{n}''(L) = 0,\ \ \ \ \ n=1,2,\ldots,\\
\endaligned
\eqno (4.7)
$$
with the initial functions $y_{0}=\alpha_{0}:=\alpha$ and $y_{0}=\beta_{0}:=\beta$ respectively satisfy
$$
\aligned
&\alpha=\alpha_{0}\leq\alpha_{1}\leq\alpha_{2}\leq\cdot\cdot\cdot\leq\alpha_{n-1} \leq \alpha_{n}\leq\cdot\cdot\cdot\leq \beta_{n}\leq \beta_{n-1}\leq\cdot\cdot\cdot\leq \beta_{2}\leq \beta_{1}\leq \beta_{0}=\beta,\\
&\alpha''=\alpha_{0}''\geq\alpha_{1}''\geq\alpha_{2}''\geq\cdot\cdot\cdot\geq\alpha_{n-1}'' \geq \alpha_{n}''\geq\cdot\cdot\cdot\geq \beta_{n}''\geq \beta_{n-1}''\geq\cdot\cdot\cdot\geq \beta_{2}''\geq \beta_{1}''\geq \beta_{0}''=\beta''\\
\endaligned
\eqno (4.8)
$$
and converge uniformly to the extremal solutions of BVP (4.1), (1.2) in $[\alpha,\beta]$.

\noindent{\bf Proof.} Define
the mapping $F: W^{4,1}[0,L]\rightarrow L^{1}[0,L]$  by
$$
F(\sigma)(x)=f(x,\int_0^L \sigma(x)dx, \sigma''(x))-M\sigma''(x)+N\int_0^L\sigma(t)dt. \eqno (4.9)
$$ Denote $\Phi=T\circ F$, where $T=G^{-1}: L^{1}[0,L]\rightarrow D(G)$.  By Section 2, it is easy to see that $T: L^{1}[0,L]\rightarrow D(G)\hookrightarrow C^{2}[0,L]\cap D(G)$ is compact, then $\Phi: C^{2}[0,L]\rightarrow C^{2}[0,L]$ is completely continuous. Obviously, the solutions of (4.1), (1.2) in $W^{4,1}[0,L]$ is equivalent to the fixed-points of the mapping $\Phi$.

Firstly, we show the following result of the method of lower and upper solutions:
 $$\alpha\leq u\leq\beta \ \ \ \text{and}\ \ \ \alpha'' \geq u''\geq \beta'' \ \Rightarrow \ \alpha\leq \Phi u\leq\beta\ \ \ \text{and}\ \ \ \alpha''\geq (\Phi u)''\geq\beta''.\eqno (4.10)$$
 Let $g=\Phi u-\alpha$, by the the lower solution given in Definition 4.1 and (4.4), (4.5), (4.6), we have
$$
\aligned
&g^{(4)}-Mg''+N\int_0^Lg(t)dt=[(\Phi u)^{(4)}-M(\Phi u)''+N\int_0^L(\Phi u)(t)dt]-[\alpha^{(4)}-M\alpha''+N\int_0^L\alpha(t)dt]\\
&=F(u)(x)-[\alpha^{(4)}-M\alpha''+N\int_0^L\alpha(t)dt] \\
&=[f(x,\int_0^Lu(t)dt, u'')-Mu''+N\int_0^Lu(t)dt]- [\alpha^{(4)}-M\alpha''+N\int_0^L\alpha(t)dt]\\
&\geq f(x,\int_0^Lu(t)dt, u'')-f(x,\int_0^L\alpha(t)dt, \alpha'')-M(u''-\alpha'')+N\int_0^L(u-\alpha)(t)dt\\
&=f(x,\int_0^Lu(t)dt, u'')-f(x,\int_0^Lu(t)dt, \alpha'')+f(x,\int_0^Lu(t)dt, \alpha'')-f(x,\int_0^L\alpha(t)dt, \alpha'')\\
&\ \ \ \ \ \ -M(u''-\alpha'')+N\int_0^L(u-\alpha)(t)dt\\
&=[M-(a+b\int_0^L u(x)dx)](\alpha''-u'')+[N-(c-b\alpha'')]\int_0^L (u-\alpha) (x)dx\\
&\geq [M-(a+b\int_0^L \beta(x)dx)](\alpha''-u'')+[N-(c-b\beta'')]\int_0^L (u-\alpha) (x)dx.\\
&\geq 0.\\
 \endaligned\eqno (4.11)
  $$
On the other hand, since
$$
g(0)=(\Phi u)(0)-\alpha(0)= 0, \ \ \ \ \ g(L)=(\Phi u)(L)-\alpha(L)=0, \eqno (4.12)
$$
$$
g''(0)=(\Phi u)''(0)-\alpha''(0)= 0, \ \ \ \ \ g''(L)=(\Phi u)''(L)-\alpha''(L)=0.\eqno (4.13)
$$ By maximum principle of operator $G$ given in Theorem 3.1,  (4.11)-(4.13) imply that $g\geq0$ and $g''\leq0$, then we conclude that $\alpha\leq \Phi u$ and $\alpha''\geq(\Phi u)''$.

 By a similar way, using the definition of the upper solutions and the maximum principle in Theorem 3.1, we can get that $\Phi u\leq\beta$ and $(\Phi u)''\geq\beta''$, then (4.10) is proved.

Based upon Schauder fixed-point theorem, $\Phi$ has at least one fixed point in $[\alpha,\beta]=\{u\in C^2[0,1]| \alpha\leq u\leq\beta,\ \ \alpha''\geq u''\geq\beta''\}$ which is a solution of (4.1), (1.2).

Now, we develop the monotone iterative technique and we first show the following claim:  $$\beta\geq u_{1}\geq u_{2}\geq\alpha\ \ \ \text{and}\ \ \ \alpha'' \geq u_{2}''\geq u_{1}''\geq \beta'' \Rightarrow \Phi u_{1}\geq\Phi u_{2}\ \ \ \text{and}\ \ \ (\Phi u_{1})''\leq (\Phi u_{2})''.\eqno (4.14)$$
In fact, let $\phi=\Phi u_{1}-\Phi u_{2}$, we have
$$
\aligned
&\phi^{(4)}-M\phi''+N\int_0^L\phi(t)dt\\
&=[(\Phi u_{1})^{(4)}-M\Phi u_{1}''+N\int_0^L(\Phi u_{1})(t)dt]-[(\Phi u_{2})^{(4)}-M\Phi u_{2}''+N\int_0^L(\Phi u_{2})(t)dt] \\
&=F(u_{1})(x)-F(u_{2})(x)\\
&=f(x,\int_0^Lu_{1}(t)dt, u_{1}'')-f(x,\int_0^Lu_{2}(t)dt, u_{2}'')-M(u_{1}''-u_{2}'')+N\int_0^L(u_{1}-u_{2})(t)dt\\
&=f(x,\int_0^Lu_{1}(t)dt, u_{1}'')-f(x,\int_0^Lu_{1}(t)dt, u_{2}'')+f(x,\int_0^Lu_{1}(t)dt, u_{2}'')-f(x,\int_0^Lu_{2}(t)dt, u_{2}'')\\
&\ \ \ \ \ \ -M(u_{1}''-u_{2}'')+N\int_0^L(u_{1}-u_{2})(t)dt\\
&=[M-(a+b\int_0^L u_{1}(x)dx)](u_{2}''-u_{1}'')+[N-(c-bu_{2}'')]\int_0^L (u_{1}-u_{2}) (x)dx\\
&\geq [M-(a+b\int_0^L \beta(x)dx)](u_{2}''-u_{1}'')+[N-(c-b\beta'')]\int_0^L (u_{1}-u_{2}) (x)dx.\\
&\geq 0.\\
 \endaligned\eqno (4.15)
  $$
On the other hand, $$
\phi(0)=\phi(L)=\phi''(0)=\phi''(L)=0,
$$
then by Theorem 3.1, we conclude that $\phi\geq0$ and $\phi''\leq0$, and the claim (4.14) is proved.

By the definition of the mapping $\Phi$, the iterative procedure (4.7) is equivalent to the iterative equation
$$
y_{n}=\Phi y_{n-1},\ \ \ \ n=1,2,\ldots.\eqno (4.16)
$$
Define the iterative sequences $\{\alpha_{n}\}$ and $\{\beta_{n}\}$ satisfying $$
\alpha_{n}=\Phi\alpha_{n-1},\ \ \ \beta_{n}=\Phi\beta_{n-1},\ \ \ \ \ n=1,2,\ldots \eqno (4.17)
$$
with $\alpha_{0}=\alpha$ and $\beta_{0}=\beta$. Then combining (4.10) with (4.14), it is easy to see that $\{\alpha_{n}\}$ and $\{\beta_{n}\}$ have the monotonicity (4.8). By the compactness of $\Phi$ and the monotonicity (4.8), it follows that $\{\alpha_{n}\}$ and $\{\beta_{n}\}$ are convergent in $C[0,L]$, that is, there exist $\underline{y}$ and $\overline{y}\in C[0,L]$ such that $$\lim\limits_{n\rightarrow \infty}\alpha_{n}(x)=\underline{y}(x), \ \ \ \ \  \ \ \ \ \lim\limits_{n\rightarrow \infty}\beta_{n}(x)=\overline{y}(x).\eqno (4.18)$$
On the other hand, it is easy to see that the operator $\Phi$ is continuous, then letting $n\rightarrow\infty$ in (4.17), we have
$$
\underline{y}=\Phi(\underline{y}),\ \ \ \ \ \ \ \ \  \overline{y}=\Phi(\overline{y}),
$$
thus $\underline{y}$ and $\overline{y}$ are the solutions of (4.1), (1.2).

Finally, we show that $\underline{y}$ and $\overline{y}$ are the extremal solutions of (4.1), (1.2) on $[\alpha, \beta].$

Let $y\in [\alpha,\beta]$ be an arbitrary solution of problem (4.1), (1.2), then combining (4.10) with (4.14) we have
$$
\Phi^{n}\alpha\leq \Phi^{n}y\leq\Phi^{n}\beta,\ \ \ \ \ \ (\Phi^{n}\alpha)''\geq (\Phi^{n}y)''\geq(\Phi^{n}\beta)''
$$
that is
$$
\alpha_{n}\leq y\leq\beta_{n},\ \ \ \ \ \ \ \ \ \alpha_{n}''\geq y\geq\beta_{n}''.
$$
Letting $n\rightarrow\infty$, we have
$$
\underline{y}\leq y\leq\overline{y},\ \ \ \ \ \ \ \ \ \underline{y}''\geq y''\geq\overline{y}''.\eqno (4.19)
$$
Hence, $\underline{y}$ and $\overline{y}$ are minimum and maximum solutions of (4.1), (1.2) in $[\alpha, \beta]$ respectively.\hfill{$\Box$}

\noindent{\bf Remark 4.1} Substituting the operator $G$
defined in (3.1) and the function $f$
 defined in (4.4) into the iterative scheme (4.7) of Theorem 4.1, we obtain its explicit form as the following algorithm sequence:
\begin{small}
$$
\aligned
&y_{n}^{(4)}-My_{n}''+N\int_0^Ly_{n}(t)dt\\
&\ \ \ \ \ \ \ \ \ \ \ \ \ =p(x)+ay_{n-1}''+by_{n-1}''\int_0^L y_{n-1}(x)dx-c\int_0^L y_{n-1}(x)dx-My_{n-1}''(x)+N\int_0^Ly_{n-1}(t)dt,\\
&y_{n}(0) = y_{n}(L) = y_{n}''(0) = y_{n}''(L) = 0,\ \ \ \ \ \ \ \ \ \ \ \ \ \ \ \ \ \ \ \ n=1,2,\ldots.\\
\endaligned
$$
\end{small}
In the above iterative algorithm, $y_{n}$ can be given explicitly by (2.9) and (2.7) in Theorem 2.1.

\noindent{\bf Remark 4.2} In the proof of Theorem 4.1, (4.10) provides an a priori global bound for the iterative sequences; (4.14) shows that the continuous iterative operator $\Phi$ given in (4.16) is order-preserving and monotone; (4.8) in the conclusion of Theorem 4.1 demonstrates that the resulting iterative sequences $\{\alpha_{n}\}$ and $\{\beta_{n}\}$ are monotonically increasing and decreasing, respectively, thereby forming a sequence of contracting intervals $[\alpha_{n},\beta_{n}]\subseteq[\alpha,\beta]$.



In Theorem 4.1, adding some conditions we can obtain following uniqueness result:

\noindent{\bf Theorem 4.2.}  In Theorem 4.1, assume that
$$a+b\int_0^L \alpha(x)dx+\frac{L^3}{4}c>\frac{L^3N}{4}+M-\frac{4}{L^2},
\eqno (4.20)$$
then the iterative sequences $\{\alpha_{n}\}$ and $\{\beta_{n}\}$ converge uniformly to the unique solution of BVP (4.1), (1.2) in $[\alpha,\beta]=\{u\in C^2[0,1]| \alpha\leq u\leq\beta,\ \alpha''\geq u''\geq\beta''\}$. Moreover, for any initial function $y_{0}\in [\alpha,\beta]$, the iterative sequence $\{y_{n}\}$ produced by (4.7) converges uniformly to this unique solution.

\noindent{\bf Proof.}
 By the proof of Theorem 4.1, BVP (4.1), (1.2) has a minimum solution $\underline{y}$ and a maximum solution  $\overline{y}$ in $[\alpha,\beta]$,
moreover, every solution of BVP (4.1), (1.2) in $[\alpha,\beta]$ satisfies (4.19). We first show $\underline{y}=\overline{y}$. For the iterative sequences $\{\alpha_{n}\}$ and $\{\beta_{n}\}$, by (4.7), (4.9) and (3.4) we have
\begin{small}
$$
\aligned
&\alpha_{n+1}''-\beta_{n+1}''=(\Phi\alpha_{n})''-(\Phi\beta_{n})''\\
&=-\int^{L}_{0}G_{M}(x,s)[F(\alpha_{n})-F(\beta_{n})](s)ds+
\frac{N\int^{L}_{0}\int^{L}_{0}G(x,s)[F(\alpha_{n})-F(\beta_{n})](s)dsdx}{1+N\int^{L}_{0}\int^{L}_{0}G(x,s)dsdx}\int^{L}_{0}G_{M}(x,s)ds \\
 \endaligned\eqno (4.21)
  $$
\end{small}
Since $\beta_{n}''-\alpha_{n}''\leq 0$ and $\alpha_{n}-\beta_{n}\leq0$, by (4.6) we have
\begin{small}
$$
\aligned
&[F(\alpha_{n})-F(\beta_{n})](s)\\
&=[f(s,\int_0^L\alpha_{n}(t)dt, \alpha_{n}'')-M\alpha_{n}''+N\int_0^L\alpha_{n}(t)dt] -[f(s,\int_0^L\beta_{n}(t)dt, \beta_{n}'')-M\beta_{n}''+N\int_0^L\beta_{n}(t)dt]\\
&=f(s,\int_0^L\alpha_{n}(t)dt, \alpha_{n}'')-f(s,\int_0^L\alpha_{n}(t)dt, \beta_{n}'')+f(s,\int_0^L\alpha_{n}(t)dt, \beta_{n}'')-f(s,\int_0^L\beta_{n}(t)dt, \beta_{n}'')\\
&\ \ \ \ \ \ -M(\alpha_{n}''-\beta_{n}'')+N\int_0^L(\alpha_{n}-\beta_{n})(t)dt\\
&=[M-(a+b\int_0^L \alpha_{n}(t)dt)](\beta_{n}''-\alpha_{n}'')+[N-(c-b\beta_{n}'')]\int_0^L (\alpha_{n}-\beta_{n}) (t)dt\\
&\leq [M-(a+b\int_0^L \beta(t)dt)](\beta_{n}''-\alpha_{n}'')+[N-(c-b\beta'')]\int_0^L (\alpha_{n}-\beta_{n}) (t)dt.\\
&\leq 0,\\
 \endaligned\eqno (4.22)
  $$
\end{small}
then according to (4.8), (4.21) and the properties of the Green functions $G_{M}$ and $G_{0}$ in (2.4)-(2.5), we have

\begin{small}
$$
\aligned
&\alpha_{n+1}''-\beta_{n+1}''\leq-\int^{L}_{0}G_{M}(x,s)[F(\alpha_{n})-F(\beta_{n})](s)ds\\
&=-\int^{L}_{0}G_{M}(x,s)\{[M-(a+b\int_0^L \alpha_{n}dt)](\beta_{n}''-\alpha_{n}'')(s)+[N-(c-b\beta_{n}''(s))]\int_0^L (\alpha_{n}-\beta_{n}) dt\}ds\\
&=\int^{L}_{0}G_{M}(x,s)\{[M-(a+b\int_0^L \alpha_{n}dt)](\alpha_{n}''-\beta_{n}'')(s)+[N-(c-b\beta_{n}''(s))]\int_0^L (\beta_{n}-\alpha_{n}) dt\}ds\\
&\leq\int^{L}_{0}G_{M}(x,s)\{[M-(a+b\int_0^L \alpha dt)](\alpha_{n}''-\beta_{n}'')(s)+[N-(c-b\alpha''(s))]\int_0^L (\beta_{n}-\alpha_{n}) dt\}ds\\
&=\int^{L}_{0}G_{M}(x,s)\{[M-(a+b\int_0^L \alpha dt)](\alpha_{n}''-\beta_{n}'')+[N-(c-b\alpha'')]\int_0^L\int_0^L G_{0}(t,\tau) (\alpha_{n}''-\beta_{n}'') d\tau dt\}ds\\
&\leq\int^{L}_{0}\frac{\tanh(\frac{\sqrt{M}L}{2})}{2\sqrt{M}}\{[M-(a+b\int_0^L \alpha dt)](\alpha_{n}''-\beta_{n}'')+[N-(c-b\alpha'')]\frac{L^{2}}{4}\int_0^L (\alpha_{n}''-\beta_{n}'') d\tau \}ds,\\
 \endaligned
 $$
 \end{small}
 thus
 $$
 \aligned
&\|\alpha_{n+1}''-\beta_{n+1}''\|_{\infty}\\
&\leq L\frac{\tanh(\frac{\sqrt{M}L}{2})}{2\sqrt{M}}[M-(a+b\int_0^L \alpha dt)+\frac{L^3}{4}(N-c+b\underset{x\in [0,L]}\max\alpha''(x))] \|\alpha_{n}''-\beta_{n}''\|_{\infty}. \\
\endaligned\eqno (4.23)
 $$
By the assumption (4.20),
 $$
 \aligned
&\rho=L\frac{\tanh(\frac{\sqrt{M}L}{2})}{2\sqrt{M}}[M-(a+b\int_0^L \alpha dt)+\frac{L^3}{4}(N-c+b\underset{x\in [0,L]}\max\alpha''(x))]\\
&\leq\frac{L \frac{\sqrt{M}L}{2}}{2\sqrt M}[M-(a+b\int_0^L \alpha dt)+\frac{L^3}{4}(N-c)]<1,\\ \endaligned\eqno (4.24)
 $$
 thus the recursive estimation (4.23) implies that
 $$
 \|\alpha_{n}''-\beta_{n}''\|_{\infty}\leq\rho^{n}\|\alpha_{0}''-\beta_{0}''\|_{\infty}\rightarrow 0, \ \ \ \ \ (n\rightarrow\infty).\eqno (4.25)
 $$
 Consequently,
 $$
 \|\alpha_{n}-\beta_{n}\|_{\infty}= \|\int_0^L G_{0}(x,s) (\alpha_{n}''(s)-\beta_{n}''(s))ds\|_{\infty}\leq\frac{L^2}{4}\|\alpha_{n}''-\beta_{n}''\|_{\infty}\rightarrow 0, \ \ \  (n\rightarrow\infty).\eqno (4.26)
 $$
Hence by (4.18), $\lim\limits_{n\rightarrow \infty}\alpha_{n}=\underline{y}= \lim\limits_{n\rightarrow \infty}\beta_{n}=\overline{y}$. Since every solution $y$ of BVP (4.1), (1.2) in $[\alpha,\beta]$ satisfies (4.19), we conclude that $\widetilde{y}:=\underline{y}=\overline{y}$ is the unique solution of BVP (4.1), (1.2) in $[\alpha,\beta]$.

Given $y_{0}\in [\alpha,\beta]$, let $\{y_{n}\}$ be the iterative sequence produced by (4.7). Then $\{y_{n}\}$ satisfies the iterative Eq. (4.16). Using (4.14) we can prove that $$\alpha_{n}\leq y_{n}\leq\beta_{n},\ \ \ \ \ \ \ \ \ \alpha_{n}''\geq y_{n}''\geq\beta_{n}'',\ \ \ \ \ \ \ \ \ \ \ \ n=0,1,2,\ldots \eqno (4.27)$$
Letting $n\rightarrow\infty$, we obtain that $$y_{n}\rightarrow \widetilde{y}.\eqno (4.28)$$\hfill{$\Box$}

\noindent{\bf Remark 4.3}  In Theorem 4.2, by (4.25) and (4.26) we can obtain the following estimate for the distance between the iterative sequences of the upper and lower solutions£º
$$
 \|\alpha_{n}-\beta_{n}\|_{\infty}\leq\frac{L^2}{4}\|\alpha_{n}''-\beta_{n}''\|_{\infty}
 \leq \frac{L^2}{4} \rho^{n}\|\alpha_{0}''-\beta_{0}''\|_{\infty}=\frac{L^2}{4} \rho^{n}\|\alpha''-\beta''\|_{\infty}\rightarrow 0, \ \ \ \ \ (n\rightarrow\infty),\eqno (4.29)
 $$
where $\rho<1$ is as defined in (4.24).

\noindent{\bf Remark 4.4} By combining (4.27) and (4.28) in Theorem 4.2 with (4.29), we can easily conclude that the iterative Algorithm defined in (4.16) is stable with respect to operator $\Phi$, i.e., it possesses T-stability (see Definition 7.1 in [59] or Definition 2.1 in [60]). This implies that the algorithm will not converge to any spurious solution that is not a fixed point, and also guarantees that any convergent sequence that asymptotically consistent with the iterative scheme must converge to the true solution.

By using a proof process similar to Theorem 4.2 and the monotonicity of iterative sequences, we can obtain the following estimate for iterative sequence $\{\alpha_{n}\}$  and iterative sequence $\{\beta_{n}\}$:

\noindent{\bf Corollary 4.1} Under the conditions of Theorem 4.1 and Theorem 4.2, the iterative sequences $\{\alpha_{n}\}$ and $\{\beta_{n}\}$ converge uniformly to the unique solution $\widetilde{y}$ of BVP (4.1), (1.2) in $[\alpha,\beta]$ and satisfy the following estimates, respectively:
$$
 \|\alpha_{n+1}-\alpha_{n}\|_{\infty}\leq\frac{L^2}{4}\|\alpha_{n+1}''-\alpha_{n}''\|_{\infty}
 \leq \frac{L^2}{4} \rho^{n}\|\alpha_{1}''-\alpha_{0}''\|_{\infty}\leq\frac{L^2}{4} \rho^{n}\|\alpha''-\beta''\|_{\infty}\rightarrow 0, \ \  (n\rightarrow\infty),\eqno (4.30)
$$
$$
 \|\alpha_{n}-\widetilde{y}\|_{\infty}\leq\frac{L^2}{4}\|\alpha_{n}''-\widetilde{y}''\|_{\infty}
 \leq\frac{L^2}{4} \frac{\rho^{n}}{1-\rho}\|\alpha''-\beta''\|_{\infty}\rightarrow 0, \ \ \ \ \ (n\rightarrow\infty),\eqno (4.31)
 $$
$$
 \|\beta_{n+1}-\beta_{n}\|_{\infty}\leq\frac{L^2}{4}\|\beta_{n+1}''-\beta_{n}''\|_{\infty}
 \leq \frac{L^2}{4} \rho^{n}\|\beta_{1}''-\beta_{0}''\|_{\infty}\leq\frac{L^2}{4} \rho^{n}\|\alpha''-\beta''\|_{\infty}\rightarrow 0, \ \  (n\rightarrow\infty),\eqno (4.32)
 $$
$$
 \|\beta_{n}-\widetilde{y}\|_{\infty}\leq\frac{L^2}{4}\|\beta_{n}''-\widetilde{y}''\|_{\infty}
 \leq\frac{L^2}{4} \frac{\rho^{n}}{1-\rho}\|\alpha''-\beta''\|_{\infty}\rightarrow 0, \ \ \ \ \ (n\rightarrow\infty),\eqno (4.33)
$$
where $\rho<1$ is as defined in (4.24).

\noindent{\bf Proof.} From the monotonicity $\alpha_{n}''-\alpha_{n-1}''\leq 0$ and $\alpha_{n-1}-\alpha_{n}\leq0$ given by (4.8), we can easily obtain (4.30) by employing an argument similar to (4.21)-(4.26) in Theorem 4.2. Now, from the inequality $$\|\alpha_{n}''-\widetilde{y}''\|_{\infty}\leq \|\alpha_{n}''-\alpha_{n+1}''\|_{\infty}+\|\alpha_{n+1}''-\widetilde{y}''\|_{\infty}
=\|\alpha_{n}''-\alpha_{n+1}''\|_{\infty}+\|(\Phi\alpha_{n})''-(\Phi\widetilde{y})''\|_{\infty}, \eqno (4.34)$$
and the estimate (derived similarly to (4.21)-(4.24))
$$
\|(\Phi\alpha_{n})''-(\Phi\widetilde{y})''\|_{\infty}\leq \rho \|\alpha_{n}''-\widetilde{y}''\|_{\infty}, \eqno (4.35)
$$
we can conclude that
$$
\|\alpha_{n}''-\widetilde{y}''\|_{\infty}\leq \frac{1}{1-\rho}\|\alpha_{n}''-\alpha_{n+1}''\|_{\infty}. \eqno (4.36)
$$
With this result and the latter part of (4.30) we have
$$
\|\alpha_{n}-\widetilde{y}\|_{\infty}
\leq\frac{L^2}{4}\|\alpha_{n}''-\widetilde{y}''\|_{\infty}\leq \frac{1}{1-\rho}\frac{L^2}{4}\|\alpha_{n}''-\alpha_{n+1}''\|_{\infty}\leq\frac{L^2}{4} \frac{\rho^{n}}{1-\rho}\|\alpha''-\beta''\|_{\infty},
$$
that is, (4.31) holds.

Analogously we can prove (4.32) and (4.33). \hfill{$\Box$}

\noindent{\bf Remark 4.5} The numerical stability of the monotone iterative algorithm is established by Corollary 4.1, Remark 4.2, and Remark 4.3. Together, they demonstrate that numerical errors generated by the iterative computations are constrained by a contracting a priori error bound whose monotonic decay inherently suppresses the growth of errors. This framework is robust: the monotonicity and continuity of the iterative operator $\Phi$ (Remark 4.2) make a violation of the monotonic order relation (4.8) by numerical perturbations highly improbable under standard computational precision. 
 As a result, the algorithm confines and damps errors within the nested intervals, thereby providing a direct guarantee of numerical stability.


\section{The typical application of the monotone  iterative technique}

According to the main result Theorem 4.1 and Definition 4.1, we can usually take $\alpha\equiv0$ as the lower solution for the BVP (4.1), (1.2). For the upper solution, when the live load $p=p(x)$ meets certain conditions, we can choose different types of nonnegative functions that satisfy the boundary conditions (1.2), such as $\beta=\lambda\sin(\frac{\pi}{L} x)$, $\beta=\lambda x(x^3-2Lx^2+L^3)$, $\beta=\lambda[(L^2-6)x+6L\frac{\sinh(x)}{\sinh(L)}-x^3]$ and $\beta=\lambda[ \frac{\sinh(L)}{2}(Lx-x^2)-(\sinh(L)-\sinh(L-x)-\sinh(x))]$  etc, where $\lambda$ is some positive constant.

 In the sequel, we take the upper solution $\beta=\lambda\sin(\frac{\pi}{L} x)$ as an typical example to illustrate the application of Theorem 4.1, other forms of upper solutions can be discussed similarly.

\noindent{\bf Theorem 5.1.} 
Let $a, b, c$ are positive constants and $p\in L^{1}[0,L]$ is nonnegative.
Assume that there exists some positive constant $\lambda$ such that
 $$p(x)\leq [\lambda\frac{\pi^{4}}{L^4}+a\lambda\frac{\pi^{2}}{L^2}+2b\lambda^2\frac{\pi}{L}]\sin(\frac{\pi}{L}x)+2\lambda c\frac{L}{\pi}, \ \ \ \text{a.e.}\  x\in (0,L),\eqno (5.1)$$
 and $$c+b\lambda\frac{\pi^{2}}{L^2}\leq \frac{4(a+2b\lambda\frac{L}{\pi})}{L^3}\sigma(\sqrt{a+2b\lambda\frac{L}{\pi}}L)\eqno (5.2)$$
where the function $\sigma$ is defined as in (3.3). Then the 
classical Melan problem (4.1), (1.2) have a lower solution $\alpha\equiv0$ and an upper solution $\beta=\lambda\sin(\frac{\pi}{L} x)$ as defined in Definition 4.1. And consequently, take $$M=a+b\int_0^L \beta(x)dx=a+2b\lambda\frac{L}{\pi},\ \ \ \  N=c-b\underset{x\in [0,L]}\min\beta''(x)=c+b\lambda\frac{\pi^{2}}{L^2},\eqno (5.3)$$ then the iterative sequences $\{\alpha_{n}\}$ and $\{\beta_{n}\}$ produced by the iterative procedure (4.7) with the initial functions $y_{0}=\alpha\equiv0$ and $y_{0}=\beta=\lambda\sin(\frac{\pi}{L} x)$ will respectively converge to the minimum solution $\underline{y}$ and the maximum solutions $\overline{y}$ of (4.1), (1.2) in $[0, \lambda\sin(\frac{\pi}{L} x)]$.

\noindent{\bf Proof.} To complete the proof, we just need to verify that all the conditions of the theorem 4.1 are satisfied. By (5.1) and Definition 4.1, it is easy to verify that $\beta=\lambda\sin(\frac{\pi}{L} x)$ is an upper solution of (4.1), (1.2) and $\alpha\equiv0$ is a lower solution. Obviously, $\alpha,\beta$ satisfy the inequality (4.5); The value of the $M, N$ taking in (5.3) show that condition (4.6) holds, and then (5.2) means that the key condition $N\leq \frac{4M}{L^3}\sigma(\sqrt{M}L)$ is satisfied. Then by Theorem 4.1, the conclusion of the theorem 4.1 holds.
\hfill{$\Box$}

\noindent{\bf Remark 5.1} In Theorem 5.1, if $p \not\equiv 0$, then $y\equiv 0$ is not a solution of (4.1), (1.2) and consequently we conclude that the minimum solution $\underline{y}$ and the maximum solutions $\overline{y}$ are both positive.

\noindent{\bf Remark 5.2} In Theorem 5.1, according to (4.7), (4.4) and (5.3), the specific iterative procedure is as follows:
\begin{small}
$$
\aligned
&Ky_{n}(x)=y_{n}^{(4)}-My_{n}''+N\int_0^Ly_{n}(t)dt= y_{n}^{(4)}-(a+2b\lambda\frac{L}{\pi})y_{n}''+(c+b\lambda\frac{\pi^{2}}{L^2})\int_0^Ly_{n}(t)dt\\
&=f(x,\int_0^L y_{n-1}(x)dx, y_{n-1}''(x))-My_{n-1}''(x)+N\int_0^Ly_{n-1}(t)dt\\
&=p(x)+ay_{n-1}''+by_{n-1}''\int_0^L y_{n-1}(x)dx-c\int_0^L y_{n-1}(x)dx-My_{n-1}''(x)+N\int_0^Ly_{n-1}(t)dt\\
&=p(x)+by_{n-1}''\int_0^L y_{n-1}(x)dx-2b\lambda\frac{L}{\pi}y_{n-1}''+b\lambda\frac{\pi^{2}}{L^2}\int_0^Ly_{n-1}(t)dt,\\
 \endaligned\eqno (5.4)
  $$
\end{small}
that is
\begin{small}
$$
\aligned
&y_{n}^{(4)}-(a+2b\lambda\frac{L}{\pi})y_{n}''+(c+b\lambda\frac{\pi^{2}}{L^2})\int_0^Ly_{n}(t)dt\\
&\ \ \ \ \ \ \ \ \ \ \ \ \ \ \ \ \ \ \ \ \ \ \ \ \ \ \ \ \ \ \ =p(x)+by_{n-1}''\int_0^L y_{n-1}(x)dx-2b\lambda\frac{L}{\pi}y_{n-1}''+b\lambda\frac{\pi^{2}}{L^2}\int_0^Ly_{n-1}(t)dt,\\
&y_{n}(0) = y_{n}(L) = y_{n}''(0) = y_{n}''(L) = 0,\ \ \ \ \ n=1,2,\ldots.\\
\endaligned
\eqno (5.5)
$$
\end{small}
where each $y_{n}$ can be obtained by (2.9) and (2.7) in Theorem 2.1.

Next, we give a uniqueness result for the solutions obtained by Theorem 5.1 based on Theorem 4.2.

\noindent{\bf Theorem 5.2.} In Theorem 5.1, assuming that $p \not\equiv 0$ and
$$\lambda\frac{L^3}{4}[\frac{2b}{\pi}+\frac{b\pi^{2}}{4}]\leq1,\ \ \ \eqno (5.6)$$ 
then the positive solution of problem (4.1), (1.2) in $[0, \lambda\sin(\frac{\pi}{L} x)]$ is unique,
and the iterative sequence $\{\alpha_{n}\}$ (with initial function $\alpha_{0}\equiv0$), $\{\beta_{n}\}$ (with initial function $\beta_{0}=\lambda\sin(\frac{\pi}{L}x$)) and $\{u_{n}\}$ (with any initial function $u_{0}\in [0, \lambda\sin(\frac{\pi}{L} x)]$) produced by (5.5) will all converge uniformly to the unique positive solution of (4.1), (1.2) in $[0, \lambda\sin(\frac{\pi}{L} x)]$. 

\noindent{\bf Proof.} 
It is easy to verify that (5.6) is equivalent to (4.20) with $\alpha\equiv0$ and $M=a+2b\lambda\frac{L}{\pi}, N=c+b\lambda\frac{\pi^{2}}{L^2}$. So, by Theorem 4.2, the conclusion of the Theorem 5.2 holds.   
\hfill{$\Box$}

We present a simple example to illustrate the application of
Theorem 5.1 and 5.2.

\noindent{\bf Example 5.1. (Table 5 in [38])} Taking the parameters $a=\frac{1}{10}, b=\frac{1}{10}, c=\frac{1}{10}$ and $L=2, p(x)\equiv\frac{1}{10}$ in (4.1), (1.2), we have
$$
\aligned
&w''''(x)-(\frac{1}{10}+\frac{1}{10}\int_0^2 w(x)dx)w''(x)+\frac{1}{10}\int_0^2 w(x)dx=\frac{1}{10},\ \ \ \ x\in (0,2),\\
&w(0)=w(2)=w''(0)=w''(2)=0.\\
\endaligned \eqno (5.7)
$$
By (5.1) we can take  $\lambda=\frac{\pi}{4}$ and then it is easy to verify that problem (5.7) have a lower solution $\alpha\equiv0$ and an upper solution $\beta=\frac{\pi}{4}\sin(\frac{\pi}{2} x)$ which satisfy (4.5) in Theorem 4.1. As in (5.3), we can take $$M=a+b\int_0^L \beta(x)dx=a+2b\lambda\frac{L}{\pi}=\frac{1}{5},\ \ \ \  N=c-b\underset{x\in [0,L]}\min\beta''(x)=c+b\lambda\frac{\pi^{2}}{L^2}=\frac{16+\pi^3}{160},$$
then by (3.3) we can calculate that $\sigma(\sqrt{M}L)=\sigma(\frac{2\sqrt{5}}{5})\approx 24.219$, and consequently, we can easily verify that
$$
N=\frac{16+\pi^3}{160}\approx0.293\leq \frac{4M}{L^3}\sigma(\sqrt{M}L)\approx2.422,
$$
that is, the condition (5.2) in Theorem 5.1 (that is, $N\leq \frac{4M}{L^3}\sigma(\sqrt{M}L)$  in Theorem 4.1) is satisfied. Moreover, the uniqueness condition (5.6) for the solution also holds. Then substituting the lower solution $y_{0}=\alpha\equiv0$ and the upper solution $y_{0}=\beta=\frac{\pi}{4}\sin(\frac{\pi}{2} x)$ of problem (5.7) into
right side of the iterative procedure (5.5) respectively, by using (2.9) and (2.7) in Theorem 2.1 we can calculate that
$$
\aligned
&\alpha_{1}\approx0.0932(5x-2.5x^2-\frac{\theta(x)}{0.04074}),\\
&=0.0932(5x-2.5x^2-\frac{\sinh (\sqrt \frac{1}{5}2)-\sinh (\sqrt \frac{1}{5}(2-x))-\sinh (\sqrt \frac{1}{5}x)}{0.04074}),\\
\endaligned \eqno (5.8)
$$
with
$\underset{x\in (0,2)}\max \alpha_{1}(x)=\alpha_{1}(1)=0.01795;$
$$
\aligned
&\beta_{1}\approx0.2739(5x-2.5x^2-\frac{\theta(x)}{0.04074}),\\
&=0.2739(5x-2.5x^2-\frac{\sinh (\sqrt \frac{1}{5}2)-\sinh (\sqrt \frac{1}{5}(2-x))-\sinh (\sqrt \frac{1}{5}x)}{0.04074}).\\
\endaligned \eqno (5.9)
$$
with $\underset{x\in (0,2)}\max \beta_{1}(x)=\beta_{1}(1)=0.05274.$
For the second iteration, substituting $\alpha_{1}$ and $\beta_{1}$  into
right side of the iterative procedure (5.5) respectively 
and consequently by using the Mean Value Theorem for integral we can get the approximation result of the second iteration are:
$$
\aligned
\alpha_{2}&\approx0.10002(5x-2.5x^2-\frac{\theta(x)}{0.04074}),\\
&=0.10002(5x-2.5x^2-\frac{\sinh (\sqrt \frac{1}{5}2)-\sinh (\sqrt \frac{1}{5}(2-x))-\sinh (\sqrt \frac{1}{5}x)}{0.04074})\\
\endaligned \eqno (5.10)
$$
with
$\underset{x\in (0,2)}\max \alpha_{2}(x)=\alpha_{2}(1)=0.019259;$
$$
\aligned
\beta_{2}&\approx0.11279(5x-2.5x^2-\frac{\theta(x)}{0.04074}),\\
&=0.11279(5x-2.5x^2-\frac{\sinh (\sqrt \frac{1}{5}2)-\sinh (\sqrt \frac{1}{5}(2-x))-\sinh (\sqrt \frac{1}{5}x)}{0.04074})\\
\endaligned \eqno (5.11)
$$
with
$\underset{x\in (0,2)}\max \beta_{2}(x)=\beta_{2}(1)=0.021718.$
The results of the two iterations are shown in Figure 2 below drawn with
MATLAB:
\begin{figure}[H]
\centering
\includegraphics[width=0.9\linewidth]{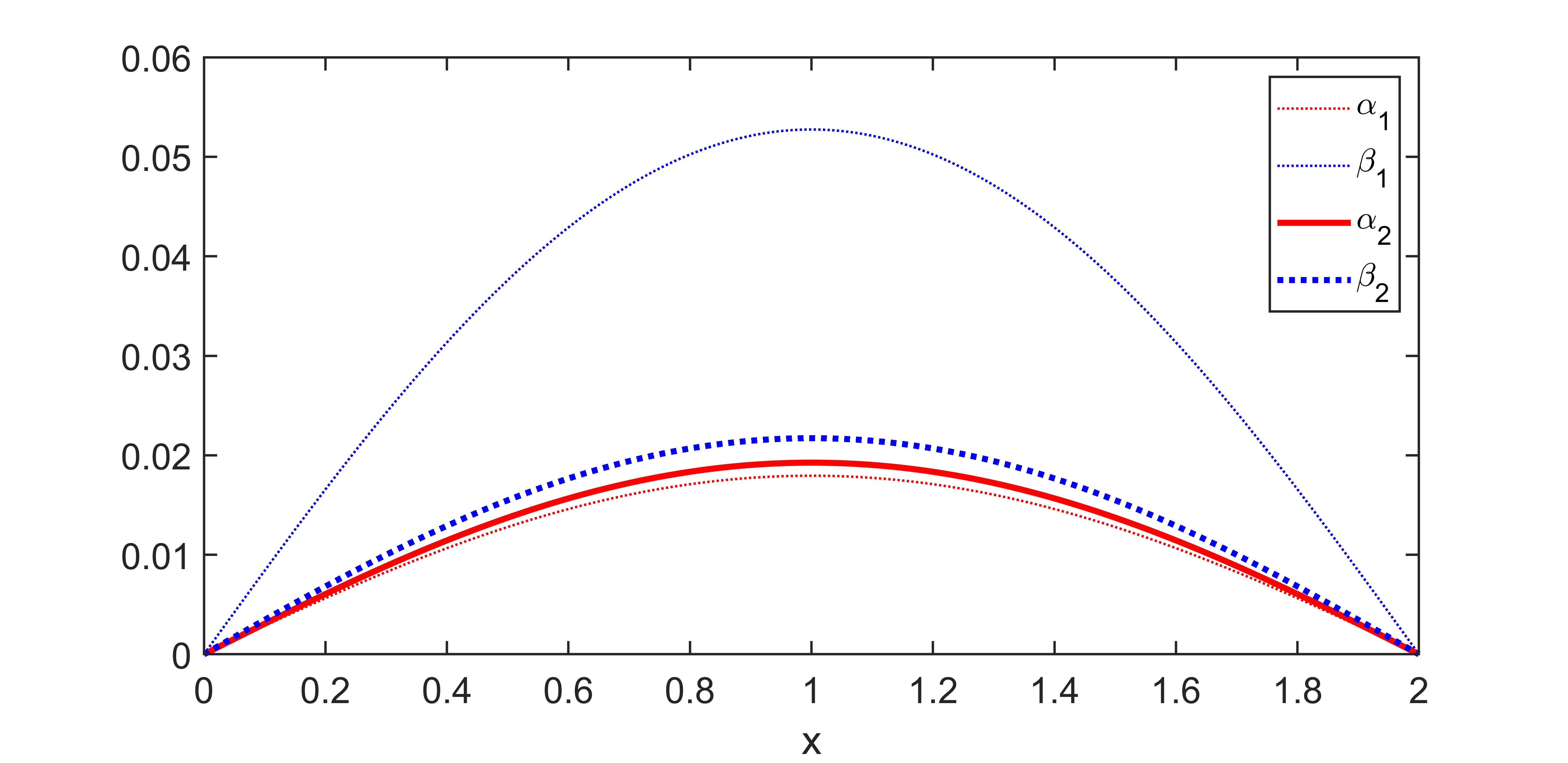}
\caption{The first and second iterations for problem (5.7).}
\end{figure}

\section{Applicability to actual bridges}

In this Section, we discuss some examples of actual bridges to demonstrate the application scope and the efficiency of the monotone iterative method obtained in Sections 2-5 in engineering design calculations. 

For the classical Melan equation (1.4) in engineering, dividing both sides by $EI$, we can obtain the equation in the form of (4.1) as follows:
$$
w''''(x)-[\frac{H}{ EI}+\frac{q}{H}\frac{E_{c}A_{c}}{EI}\frac{1}{L_{c}}\int_0^L w(x)dx]w''(x)+(\frac{q}{H})^{2}\frac{E_{c}A_{c}}{EI}\frac{1}{L_{c}}\int_0^L w(x)dx=\frac{p(x)}{EI},\ \ \ \ x\in (0,L)\eqno (6.1)
$$
that is
$$
\aligned
&a=\frac{H}{EI},\ \ \ \ \ \ b=\frac{q}{H}\frac{E_{c}A_{c}}{EI}\frac{1}{L_{c}},\ \ \ \ &c=(\frac{q}{H})^2\frac{E_{c}A_{c}}{EI}\frac{1}{L_{c}}=\frac{q}{H}b.
\endaligned \eqno (6.2)
$$
In practical engineering, the value of $\frac{q}{H}$ can be determined by the bridge design parameter $n$ called the sag-span ratio which usually ranges from $\frac{1}{5}$ to $\frac{1}{14}$ for existing bridges, and the specific relationship is: $$\frac{q}{H}=\frac{8n}{L},\eqno (6.3)$$
see e.g. [8-27]. Since the maximum possible vertical displacement of the deck of bridge is around $\frac{1}{100}$ of the length of the main span $L$ in practice (see e.g. Part 4 of Gazzola [38]), then as in Theorem 5.1, we can take $\alpha\equiv0$ as the lower solution and $$\beta=\lambda\sin(\frac{\pi}{L} x)=\frac{L}{100}\sin(\frac{\pi}{L} x)\eqno (6.4)$$ as the upper solution of (6.1). In this case, according to (5.1) and (6.3), it can be calculated that the live load $p(x)$ in equation (1.4) and (6.1) just need to satisfy:
$$
\aligned
p(x)&\leq [\frac{\pi^{4}}{100L^3}EI+\frac{H\pi^{2}}{100L}+\frac{q}{H}\frac{E_{c}A_{c}}{L_{c}}\frac{\pi L}{5000}]\sin(\frac{\pi}{L}x)+2(\frac{q}{H})^{2}\frac{E_{c}A_{c}}{L_{c}}\frac{L^{2}}{100\pi}\\
&=[\frac{\pi^{4}}{100L^3}EI+\frac{H\pi^{2}}{100L}+\frac{n\pi}{625}\frac{E_{c}A_{c}}{L_{c}}]\sin(\frac{\pi}{L}x)+\frac{32n^{2}}{25\pi}\frac{E_{c}A_{c}}{L_{c}}, \ \ \ \ \ \ \ \text{a.e.}\  x\in (0,L).\endaligned \eqno (6.5)$$
Now, we investigate the key condition (5.2) (or (3.2)) for (6.1). According to (6.3) and $\lambda=\frac{L}{100}$ in (6.4) and the relationship of $c=\frac{q}{H}b$ in (6.2), we have for actual bridges the condition (5.2) is equivalent to
$$c(\frac{800n+\pi^2}{800n})\leq \frac{4(a+\frac{cL^{3}}{400n\pi})}{L^3}\sigma(\sqrt{a+\frac{cL^{3}}{400n\pi}}L)= \frac{4(\sqrt{a+\frac{cL^{3}}{400n\pi}}L)^{2}}{L^5}\sigma(\sqrt{a+\frac{cL^{3}}{400n\pi}}L).\eqno (6.6)$$ On the other hand, according to (3.3) it is easy to verify that the function $x^{2}\sigma(x)$ is decreasing and $x^{2}\sigma(x)\leq 20, \ \forall x>0$, then (5.3), (3.2) and (5.2) imply that
$$c+b\lambda\frac{\pi^{2}}{L^2}=N\leq \frac{4M}{L^3}\sigma(\sqrt{M}L)=\frac{4}{L^5}(\sqrt{M}L)^{2}\sigma(\sqrt{M}L)\leq \frac{80}{L^5}.\eqno (6.7)$$
Substituting (6.2), (6.3) and (6.4) to the left-hand side of (6.7) we have
$$c+b\lambda\frac{\pi^{2}}{L^2}=(\frac{q}{H}+\frac{L}{100}\frac{\pi^{2}}{L^2})b=(\frac{8n}{L}+\frac{L}{100}\frac{\pi^{2}}{L^2})b\leq\frac{80}{L^5},\eqno (6.8)$$
and consequently, we can obtain
$$(\frac{800n+\pi^{2}}{100})b\leq\frac{80}{L^4}\eqno (6.9)$$
and its equivalent inequality
$$(\frac{800n+\pi^{2}}{8n})c\leq\frac{80}{L^5},\eqno (6.10)$$
which can respectively be regarded as necessary condition for the monotone iteration method we have established in Theorem 5.1 (or Theorem 4.1) to be applicable to the actual suspension bridge. On the other hand, substituting $\lambda=\frac{L}{100}$ to the uniqueness condition (5.6) in Theorem 5.2, it is easy to verify that (5.6) is equivalent to the following inequality:
$$(\frac{\frac{40}{\pi}+5\pi^{2}}{100})b\leq\frac{80}{L^4}.\eqno (6.11)$$
Comparing (6.9) with (6.11), we can come to the following conclusion: if $n\geq\frac{1}{15.32}$, then Condition (6.9) can guarantee that the uniqueness Condition (6.11) holds. Combining this fact with (6.6), we have the following uniqueness result for iterative solutions:

\noindent{\bf Theorem 6.1.} For the Melan equation (1.4), (1.2) corresponding to actual bridges, assuming that the live load $p$ satisfies (6.5), then (1.4), (1.2) have a lower solution $\alpha\equiv0$ and a upper solution $\beta=\frac{L}{100}\sin(\frac{\pi}{L} x)$. If the sag-span ratio $n\geq\frac{1}{15.32}$ and $$N=c(1+\frac{\pi^2}{800n})\leq \frac{4(\sqrt{a+\frac{cL^{3}}{400n\pi}}L)^{2}}{L^5}\sigma(\sqrt{a+\frac{cL^{3}}{400n\pi}}L)
=\frac{4M}{L^3}\sigma(\sqrt{M}L).\eqno (6.12)$$
where $a,c$ is as in (6.2), then the iterative sequences $\{\alpha_{n}\}$ and $\{\beta_{n}\}$ produced by the iterative procedure (5.5) (in which $p(x)$ should be replaced by $\frac{p(x)}{EI}$) with the initial functions $y_{0}=\alpha\equiv0$ and $y_{0}=\beta=\frac{L}{100}\sin(\frac{\pi}{L} x)$ will respectively converge to the unique positive solution of (1.4), (1.2) in $[0, \frac{L}{100}\sin(\frac{\pi}{L} x)]$.

\noindent{\bf Remark 6.1} In practical engineering, the sag-span ratio $n$ is usually designed around $\frac{1}{10}$ and ranges from $\frac{1}{5}$ to $\frac{1}{14}$, see e.g. [1], [5] or [8-27], then by Theorem 6.1, as long as the monotone iterative method is applicable, then the obtained positive solution is unique.

\noindent{\bf Remark 6.2} By the relationship of $c=\frac{q}{H}b$ in (6.2) and (6.3), the condition (6.12) is equivalent to
$$\frac{800n+\pi^2}{100}b\leq \frac{4(\sqrt{a+\frac{bL^{2}}{50\pi}}L)^{2}}{L^4}\sigma(\sqrt{a+\frac{bL^{2}}{50\pi}}L)\eqno (6.13)$$
where $a,b$ is as in (6.2).

Since the function $x^{2}\sigma(x)$ is decreasing for $x>0$, the right-hand side of inequality (6.12) decreases monotonically with respect to $c$ and $L$, respectively. Combined with the monotonic increase of the left-hand side of (6.12) with $c$, this yields the following conclusion: in order for the condition (6.12) to hold, the values of the parameter $L$ and $c$ should be relatively small, which also means that (6.12) may hold for bridges with short main span and high bridge deck stiffness $EI$. 

However, when we validate the key condition (6.12) of the monotone iterative method for actual bridges, we find that the condition appears overly restrictive. We next discuss the limitations of (6.12) and its scope of application through a practical bridge case study: 

\noindent{\bf Example 6.1.} We consider a suspension bridge across the Waal River in Nijmegen city, the Netherlands,  whose feasibility study and design schemes have been discussed in [61]. One of the reference design scheme proposed in [61] is as follows: the main span $L$ of the suspension bridge is designed to be $150$m, and the sag-span ratio is $n=\frac{1}{5}$, combining these with other material parameters (see Table 5-10 in [61]) of the bridge, we can calculate the parameters in (6.2) as follows:
$$
\aligned
&a=\frac{H}{EI}=\frac{16542}{4.557\times10^{8}}=3.63\times10^{-5},\\ &b=\frac{q}{H}\frac{E_{c}A_{c}}{EI}\frac{1}{L_{c}}=0.01066\frac{4.12173\times10^{6}}{4.557\times10^{8}}\frac{1}{163.2}=5.5\times10^{-7},\\
&c=(\frac{q}{H})^2\frac{E_{c}A_{c}}{EI}\frac{1}{L_{c}}=\frac{q}{H}b=5.5\times10^{-9}.
\endaligned \eqno (6.14)
$$
As in (6.4), we can take $\alpha\equiv0$ as the lower solution and $\beta=\lambda\sin(\frac{\pi}{L} x)=\frac{150}{100}\sin(\frac{\pi}{L} x)$ as the upper solution by assuming the live load $p=p(x)$ meets condition (6.5). We validate the key condition (6.12) or (5.2) and calculate that
$$
\aligned
&N=c+b\lambda\frac{\pi^{2}}{L^2}=\frac{800n+\pi^2}{800n}c=5.839\times10^{-9},\\  &\frac{4M}{L^3}\sigma(\sqrt{M}L)=\frac{4M}{L^3}\sigma(\sqrt{a+2b\lambda\frac{L}{\pi}}L)
=\frac{4(\sqrt{a+\frac{cL^{3}}{400n\pi}}L)^{2}}{L^5}\sigma(\sqrt{a+\frac{cL^{3}}{400n\pi}}L)\\
&\ \ \ \ \ \ \ \ \ \ \ \ \ \ \ \ \ \ \ \ \ \ \ \ \ \ \ \ \ \ \ \ \ \ \ \ \ \ \ \ \ \ \ \ \ \ \ \  =1.01\times10^{-9}.
\endaligned \eqno (6.15)$$
(6.15) shows that the key condition (6.12) or (5.2) dose not hold, that is,  our monotone iterative method in Theorem 6.1 and Theorem 5.1 appears to be inapplicable to this reference design scheme.

Next, we adjust some structural parameters of Waal River suspension bridge within practically feasible ranges to explore the applicability scope of the condition. We demonstrate this through a separate example as follows: 

\noindent{\bf Example 6.2.} 
If we design the suspension bridge in Example 6.1 to have a shorter main span of $L=100m$ (with a corresponding increase in side span length), and reduce the original value of the sag-span ratio $n$ from $\frac{1}{5}$ to $\frac{1}{9}$ (this is feasible in practical design), 
assuming other material-structural parameters are the same as before, then we have: $\frac{q}{H}=n\frac{8}{L}=0.008888$, the parameter $H=\frac{q}{0.008888}=19850.4$, $L_{c}=
\frac{L}{2}\sqrt{1+\frac{L^{2}}{4}(\frac{q}{H})^{2}}+\frac{H}{q}\ln(\frac{L}{2}\frac{q}{H}+\sqrt{1+\frac{L^{2}}{4}(\frac{q}{H})^{2}})
=103.2m$ (see (2.2) in [38]). And consequently, the following results can be obtained through calculation:
$$
\aligned
&a=\frac{H}{EI}=\frac{19850.4}{4.557\times10^{8}}=4.356\times10^{-5},\\ &b=\frac{q}{H}\frac{E_{c}A_{c}}{EI}\frac{1}{L_{c}}=0.00888\frac{4.1217\times10^{6}}{4.557\times10^{8}}\frac{1}{103.2}=7.7898\times10^{-7},\\
&c=(\frac{q}{H})^2\frac{E_{c}A_{c}}{EI}\frac{1}{L_{c}}=\frac{q}{H}b=6.9236\times10^{-9}.
\endaligned \eqno (6.16)
$$
Then the two sides of the condition (6.12) (or (5.2)) are:
$$
\aligned
&N=c+b\lambda\frac{\pi^{2}}{L^2}=7.6916\times10^{-9},\\  &\frac{4M}{L^3}\sigma(\sqrt{M}L)=\frac{4M}{L^3}\sigma(\sqrt{a+2b\lambda\frac{L}{\pi}}L)
=\frac{4\times0.000093}{100^3}\sigma(0.9653)\\
&\ \ \ \ \ \ \ \ \ \ \ =3.727\times10^{-10}\times 20.6858=7.7097\times10^{-9},
\endaligned \eqno (6.17)$$
that is, the key condition (6.12) (or (5.2)) is satisfied, and our monotone iterative method in Theorem 6.1 (or Theorem 5.1) is applicable to this design scheme. In this case, according to condition (6.5), it can be calculated that the corresponding actual live load (kN/m) just need to satisfy the following:
 $$p(x)\leq 484.8581\sin(\frac{\pi}{100}x)+223.7782, \ \ \ \text{a.e.}\  x\in (0,100),\eqno (6.18)$$
which can cover all possible live-load situations in actual traffic conditions. 
 For simplicity, we assume without loss of generality that a large actual live load $p(x)\equiv200$kN/m is uniformly distributed, we now calculate the corresponding deformation of the bridge in Example 6.2 by using the monotone iterative technique developed in Theorem 5.1, 5.2 and 6.1. In this case, the problem corresponding to (6.1) is as follows:
$$
\aligned
&w''''(x)-(a+b\int_0^{100} w(x)dx)w''(x)+c\int_0^{100} w(x)dx=\frac{200}{4.557\times10^{8}},\ \ \ \ x\in (0,100),\\
&w(0)=w(100)=w''(0)=w''(100)=0,\\
\endaligned \eqno (6.19)
$$
where $a,b,c$ is as in (6.16).
Substituting the lower solution $y_{0}=\alpha\equiv0$ and the upper solution $y_{0}=\beta=\frac{100}{100}\sin(\frac{\pi}{100} x)=\sin(\frac{\pi}{100} x)$ of problem (6.19) into
right side of the iterative procedure (5.5) respectively, by using (2.9) and (2.7) in Theorem 2.1 we can calculate that
$$
\aligned
&\alpha_{1}= 2.76769\times10^{-7}(536615.23167x-5366.1523167x^2-\frac{\theta(x)}{9.74389\times10^{-9}}),\\
&=0.1485x-0.00149x^2-28.4[\sinh (0.96488)-\sinh ( \sqrt{0.000093}(100-x))-\sinh (\sqrt{0.000093}x)]\\
\endaligned \eqno (6.20)
$$
with
$$\underset{x\in (0,100)}\max \alpha_{1}(x)=\alpha_{1}(50)=0.3342m;$$
$$
\aligned
&\beta_{1}= 3.07619\times10^{-7}(536615.23167x-5366.1523167x^2-\frac{\theta(x)}{9.74389\times10^{-9}}),\\
&=0.165x-0.00165x^2-31.57[\sinh (0.96488)-\sinh ( \sqrt{0.000093}(100-x))-\sinh (\sqrt{0.000093}x)]\\
\endaligned \eqno (6.21)
$$
with
$$\underset{x\in (0,100)}\max \beta_{1}(x)=\beta_{1}(50)=0.3709m.$$
For the second iteration, substituting $\alpha_{1}$ and $\beta_{1}$  into
right side of the iterative procedure (5.5) 
and then by the Mean Value Theorem for integral we have the results of the second iteration are
\begin{small}
$$
\aligned
&\alpha_{2}\approx2.91572\times10^{-7}(536615.23167x-5366.1523167x^2-\frac{\theta(x)}{9.74389\times10^{-9}}),\\
&=0.156462x-0.00156462x^2-29.9236[\sinh (0.96488)-\sinh ( \sqrt{0.000093}(100-x))-\sinh (\sqrt{0.000093}x)]\\
\endaligned \eqno (6.22)
$$
\end{small}
with$$\underset{x\in (0,100)}\max \alpha_{2}(x)=\alpha_{2}(50)=0.35151m;$$
\begin{small}
$$
\aligned
&\beta_{2}\approx2.9296548703\times10^{-7}(536615.23167x-5366.1523167x^2-\frac{\theta(x)}{9.74389\times10^{-9}}),\\
&=0.157209x-0.00157209x^2-30.0666[\sinh (0.96488)-\sinh ( \sqrt{0.000093}(100-x))-\sinh (\sqrt{0.000093}x)]\\
\endaligned \eqno (6.23)
$$
\end{small}
with
$$\underset{x\in (0,100)}\max \beta_{2}(x)=\beta_{2}(50)=0.35319m.$$
The results of the two iterations are shown in Figure 3 below drawn with
MATLAB:
\begin{figure}[H]
\centering
\includegraphics[width=0.9\linewidth]{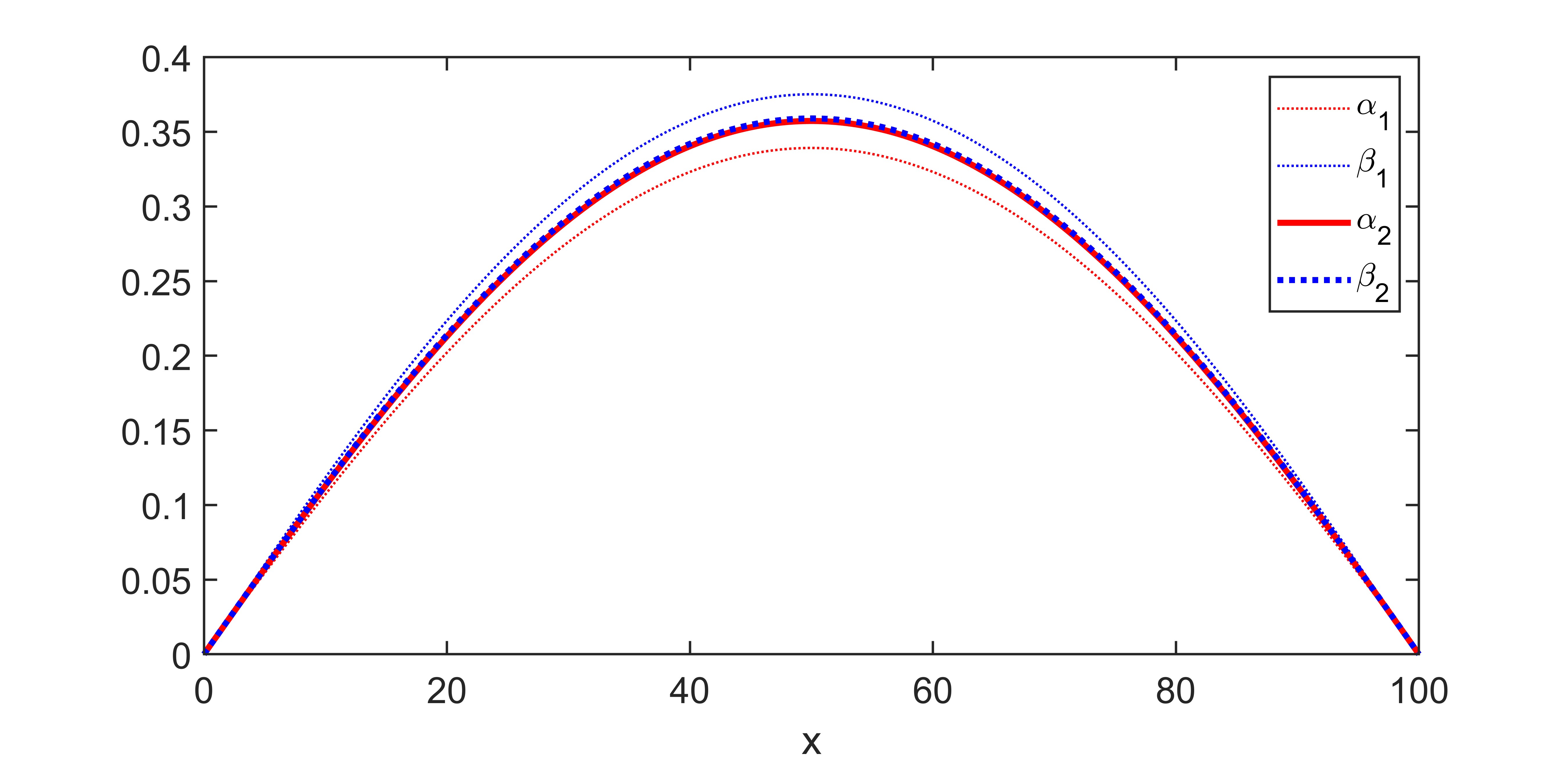}
\caption{The first and second iterations for problem (6.19).}
\end{figure}


\noindent{\bf Remark 6.3.}  In Section 5, If we take some other type of function like $\beta=\lambda x(x^3-2Lx^2+L^3)$ or $\beta=\lambda[ \frac{\sinh(L)}{2}(Lx-x^2)-(\sinh(L)-\sinh(L-x)-\sinh(x))]$ as the upper solution, we can obtain different expressions of the key condition corresponding to (5.2) and (6.12) following a similar discussion as in Section 5 and Section 6. Nonetheless, it can be verified that they still require that the values of the parameter $L$ and $c$ be relatively small, that is, the practical application scope still seems to be limited to bridges with short main span and high bridge deck stiffness $EI$. 

Next, we use an example of a real bridge to demonstrate that even if the conditions are not satisfied, our iterative sequence is still effective for  some specific load conditions.


\noindent{\bf Example 6.3.} We consider the longer actual bridge given by Wollmann [14] which has also been discussed in Gazzola [38]. The main span $L$ of the suspension bridge is $460$m, and the sag-span ratio is $f=\frac{1}{10}$, combining these with other material parameters of the bridge, we can calculate the parameters in (6.2) as follows:
$$
\aligned
&a=\frac{H}{EI}=\frac{97750}{5.7\times 10^{7}}=0.00171,\ \ b=\frac{q}{H}\frac{E_{c}A_{c}}{EI}\frac{1}{L_{c}}=0.001739\frac{3.6\times 10^{7}}{5.7\times 10^{7}}\frac{1}{472}=2.326\times 10^{-6},\\
&c=(\frac{q}{H})^2\frac{E_{c}A_{c}}{EI}\frac{1}{L_{c}}=\frac{q}{H}b=4.045\times10^{-9}.
\endaligned \eqno (6.24)
$$
By (6.4), we can take $\alpha\equiv0$ as the lower solution and $\beta=\lambda\sin(\frac{\pi}{L} x)=\frac{460}{100}\sin(\frac{\pi}{460} x)$ as the upper solution by assuming the live load $p=p(x)$(kN/m) meets the following condition as in (6.5) or (5.1):
$$
\aligned
p(x)&\leq [\frac{\pi^{4}}{100L^3}EI+\frac{H\pi^{2}}{100L}+\frac{q}{H}\frac{E_{c}A_{c}}{L_{c}}\frac{\pi L}{5000}]\sin(\frac{\pi}{L}x)+2(\frac{q}{H})^{2}\frac{E_{c}A_{c}}{L_{c}}\frac{L^{2}}{100\pi}\\
&=59.87855\sin(\frac{\pi}{460}x)+310.7102, \ \ \ \ \ \ \ \text{a.e.}\  x\in (0,460).\endaligned \eqno (6.25)$$

By (5.2), (5.3) and (6.24) we can obtain
$$
\aligned
&N=c+b\lambda\frac{\pi^{2}}{L^2}=4.535\times10^{-9},\ \ \ \ \ M=a+2b\lambda\frac{L}{\pi}=0.0048,\\  &\frac{4M}{L^3}\sigma(\sqrt{M}L)=\frac{4M}{L^3}\sigma(\sqrt{a+2b\lambda\frac{L}{\pi}}L)=\frac{4\times0.0048}{460^3}\sigma(32)\\
&\ \ \ \ \ \ \ \ \ \ \ \ \ \ \ \ \ \ =1.991\times10^{-10}\times 1.383\times10^{-11}=2.754\times10^{-21}.
\endaligned \eqno (6.26)$$
(6.26) shows that the key condition (6.12) (or (5.2)) does not hold, that is,
Theorem 6.1 and Theorem 5.1 do not seem to be applicable to this bridge under arbitrary possible practical live load satisfying (6.25). However, we find that our monotone iterative technique as presented in (5.5) still seems to be applicable and efficient for some specific loads. Considering a freight train of length $460$m having a weight density of $20$kN/m, that is $p(x)\equiv20$kN/m is uniformly distributed. In this case, the problem corresponding to (6.1) is as follows:
$$
\aligned
&w''''(x)-(a+b\int_0^{460} w(x)dx)w''(x)+c\int_0^{460} w(x)dx=\frac{20}{5.7\times 10^{7}},\ \ \ \ x\in (0,460),\\
&w(0)=w(460)=w''(0)=w''(460)=0,\\
\endaligned \eqno (6.27)
$$
where $a,b,c$ is as in (6.24).
Substituting the lower solution $y_{0}=\alpha\equiv0$ and the upper solution $y_{0}=\beta=\frac{460}{100}\sin(\frac{\pi}{460} x)=4.6\sin(\frac{\pi}{460} x)$ of problem (6.27) into
right side of the iterative procedure (5.5) respectively, by using (2.9) and (2.7) in Theorem 2.1 we can calculate that
$$
\aligned
&\alpha_{1}= 4.0901167267\times10^{-8}(47916.6666667x-104.166666667x^2-\frac{\theta(x)}{798550651.221}),\\
&=0.00195984759821x-0.00000426053825699x^2\\
&\ \ \ \ \ \ \ -5.1219252285\times10^{-17}[\sinh (31.8697348593)-\sinh ( \sqrt{0.0048}(460-x))-\sinh (\sqrt{0.0048}x)]\\
\endaligned \eqno (6.28)
$$
with
$$\underset{x\in (0,460)}\max \alpha_{1}(x)=\alpha_{1}(230)=0.223607249946m;\eqno (6.29)$$
$$
\aligned
&\beta_{1}= 4.1037448728\times10^{-8}(47916.6666667x-104.166666667x^2-\frac{\theta(x)}{798550651.221}),\\
&=0.00196637775155x-0.00000427473424251x^2\\
&\ \ \ \ \ \ \ -5.1389913295\times10^{-17}[\sinh (31.8697348593)-\sinh ( \sqrt{0.0048}(460-x))-\sinh (\sqrt{0.0048}x)]\\
\endaligned \eqno (6.30)
$$
with
$$\underset{x\in (0,460)}\max \beta_{1}(x)=\beta_{1}(230)=0.22435m.\eqno (6.31)$$

The results of the first iteration are shown in Figure 4 below drawn with
MATLAB:
\begin{figure}[H]
\centering
\includegraphics[width=0.9\linewidth]{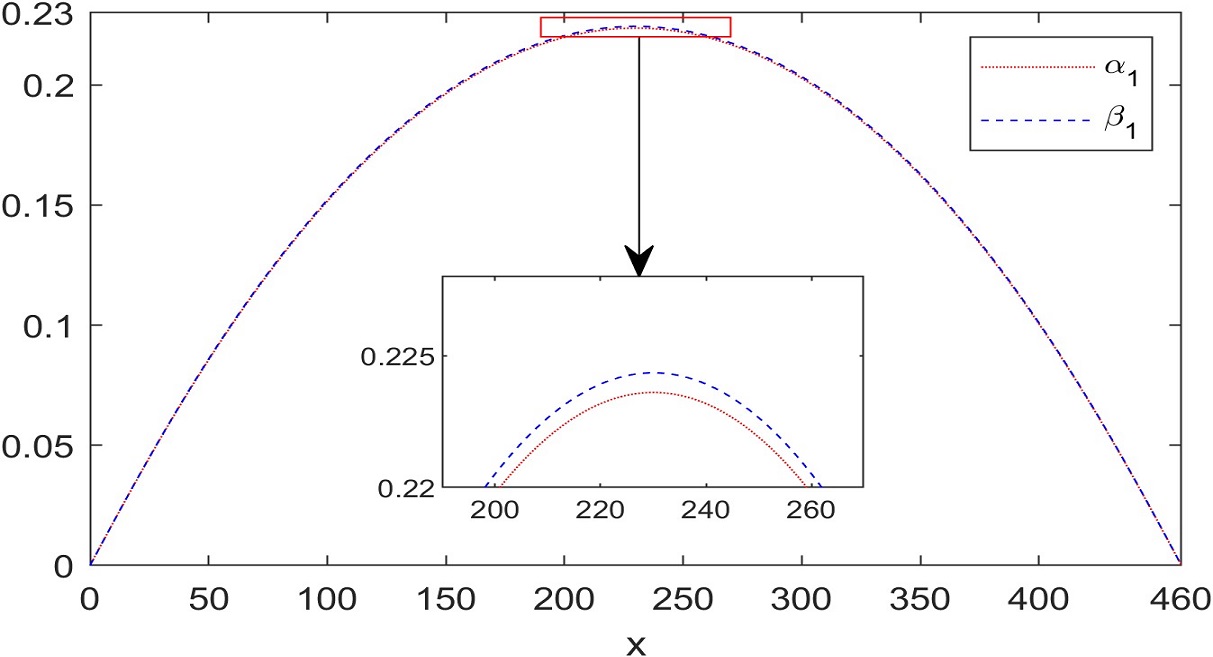}
\caption{The first iteration for problem (6.27).}
\end{figure}

\noindent{\bf Remark 6.4.} By (6.28)-(6.31) and Figure 4 in Example 6.3., it appears that our monotone iteration technique is still applicable and highly efficient for longer suspension bridges under some specific loads such as uniformly distributed load.

\section{Conclusions and open problems}

Although several approximate methods for the solution of the fundamental equations (1.1), (1.2) and (1.4), (1.2) have been proposed since the deflection theory of suspension bridges was established more than 100 years ago, they all lack rigorous proof of convergence. In this paper, based on the study of the exact explicit expression and the uniform positivity of the analytical solution of simplified ``less stiff'' model (1.5), (1.2), we develop a monotone iterative technique of upper and lower solutions to prove the existence, uniqueness and approximability of the solution for the fundamental equation (1.4), (1.2). Our approach simultaneously addresses the dual gaps, namely the absence of research on the analytical solution for the simplified ``less stiff'' model (1.5), (1.2) and the lack of a provably convergent approximation method for the original nonlinear model (1.4), (1.2). The exact explicit expression we obtained for the unique analytical solution of (1.5), (1.2) (see Theorem 2.1 and Corollary 2.1) not only provides a crucial theoretical understanding of the fundamental behavior of this simplified ``less stiff'' model but also establishes an exact benchmark, thus providing the definitive criterion for evaluating the performance of numerical methods developed in [28-33]; By analyzing the uniform positivity of the analytical solution, we established a maximum principle for the linear integro-differential operator and then constructed two successively monotone iterative sequences, which converge uniformly to the extremal solutions of the original nonlinear model (1.4), (1.2) between the lower and upper solutions. Besides, we proved that for parameters lying in some suitable
range, the two monotone sequences will converge uniformly to the unique solution of (1.4), (1.2). The monotone iterative algorithm developed in this paper provides partial answers to open problems proposed by Gazzola et al. in [38]. Our main results (see Theorem 4.1 and its typical case Theorem 5.1) showed that: the solution of the simplified ``less stiff'' model (1.5), (1.2) with appropriate modification can be used to iteratively approximate the solution of the original fundamental equation (1.4), (1.2) under certain conditions.


The applicability of the monotone iterative technique for practical engineering design is discussed in Section 6. By verifying some examples of actual bridges we showed that: for suspension bridges with short main span $L$ and high bridge deck stiffness $EI$, our approach can guarantee the uniqueness of the solution and approximate the deflection curve of the bridge under various possible practical live load situations with fast convergence rate, see Example 6.2, Figure 3 and Theorem 6.1; for existing bridges with longer main span, although the material-structural parameters usually fail to satisfy the key condition (6.12) or (5.2), our monotone iterative technique seems to be still applicable and efficient for some specific load situations such as uniformly distributed load, see Example 6.3 and Figure 4. While the key conditions (6.12) and (5.2) in Section 5 and Section 6 are only special representations of (3.2) in our method when taking special form of upper solution $\beta=\lambda\sin(\frac{\pi}{L} x)$, it seems that we cannot significantly relax the restriction of the key condition in engineering applications by replacing other forms of the upper solution functions, see Remark 6.3. This gives rise to two natural questions. Can the limited applicability of our method in practical engineering be extended? If the limitation is essential, 
what is the maximum applicable scope of our iterative method in practical engineering when the key condition is satisfied or violated, respectively?

We are confident that this paper might be the starting point for analytical approximation of the Melan equation, and we believe that some further research is needed to fully understand the role of nonlinear nonlocal term and the relationship between the simplified ``less stiff'' model and the original classical nonlinear nonlocal model.

\section{ Appendix: The derivation of classical Melan equation}

Following Von Karman and Biot [21, Section VII.1], the deflection theory of suspension bridges models the bridge structure as a combined
system of a flexible string (the sustaining cable) and a beam (the deck), see Figure 1 in Section 1. 
The beam and the string are connected through a large number of inextensible hangers. The point $O$ is the origin of the orthogonal coordinate system and positive displacements are defined as downwards. The point $M$ has coordinates $(0,L)$  with $L$ being the distance between the two towers. If no live loads act on the beam, the cable is in position $y(x),$ while the unloaded beam is in the horizontal position of the segment connecting $O$ and $M$. 
The cable is considered as a perfectly flexible string under vertical dead and live loads. When the string is subjected to a downwards vertical dead load $q,$ the horizontal component $H>0$ of the tension force remains constant. If the mass of the cable is neglected, then the dead load is distributed per horizontal unit. Hence, there is an equilibrium position in the system and the configuration of the cable is obtained by solving the equation (see [21, (1.3), Section VII])
$$
Hy''(x)=-q.\eqno (8.1)
$$
According to the solution of (8.1), the length $L_{c}$ of the cable with parabolic shape is given by
$$
L_{c}=\int_0^L\sqrt{1+(y'(x))^{2}}dx.\eqno (8.2)
$$
The function $w$ denotes both the downwards displacements of the beam and the cable, since the elastic deformation of the hangers is neglected. When the live load $p$ is imposed, a certain amount $p_{1}$ of $p$ is carried by the cable whereas the residual component $p-p_{1}$ is resisted by the bending stiffness of the beam. Then, the governing equation for the beam's downward displacement $w$ is classically expressed as:
$$
EI w''''(x)=p(x)-p_{1}(x),\ \ \ \ x\in (0,L),
\eqno (8.3)
$$
At the same time, the horizontal tension of the cable is increased to $H+h(w)$ and the deflection $w$ is coupled with the primary displacement $y$. Hence, in accordance with (8.1), the governing equation accounting for this condition reads
$$
(H+h(w))(y''(x)+w''(x))=-q-p_{1}(x),\ \ \ \ x\in (0,L),
\eqno (8.4)
$$
Then, combining (8.1), (8.3) and (8.4) we have
$$
EI w''''(x)-(H+h(w))w''(x)+\frac{q}{H}h(w)=p(x),\ \ \ \ x\in (0,L),
\eqno (8.5)
$$
which is widely recognized in literature as the Melan equation [4, p.77].  The beam representing the bridge is considered to be hinged at its endpoints, which means that the boundary conditions corresponding to (8.5) are
$$
w(0)=w(L)=w''(0)=w''(L)=0.\eqno (8.6)
$$
The elongation $\triangle L_{c}$ of the length $L_{c}$ of the cable caused by the deformation $w$ is
$$
\triangle L_{c}=\Gamma(w):=\int_0^L(\sqrt{1+(y'(x)+w'(x))^{2}}-\sqrt{1+(y'(x))^{2}})dx.\eqno (8.7)
$$
In literature, there are at least three different ways to approximate $\Gamma(w)$. Here we focus on the the way proposed by Von Karman and Biot [21, (5.14)] given its widespread adoption in engineering literature and textbooks.

Refer to the asymptotic expansion, for any $\rho\neq 0$ the following holds
$$
\sqrt{1+(\rho+\varepsilon)^{2}}-\sqrt{1+(\rho)^{2}}\sim \frac{\varepsilon\rho}{\sqrt{1+(\rho)^{2}}},  \ \ \ \text{as} \ \varepsilon\rightarrow 0.
$$
Then by applying it to (8.7) one obtains
$$
\triangle L_{c}\approx \int_0^L \frac{w'(x)y'(x)}{\sqrt{1+(y'(x))^{2}}}dx.\eqno (8.8)
$$
Following (8.8), von Karman and Biot [21, (5.14)] neglected $(y'(x))^{2}$ in comparison with $1$ and formulated
$$
\triangle L_{c}\approx \int_0^L w'(x)y'(x)dx=-\int_0^L w(x)y''(x)dx=\frac{q}{H}\int_0^L w(x)dx,\eqno (8.9)
$$
in which the integration by parts accounting for the boundary conditions $w(0)=w(L)=0$, and (8) confirming the second equality.

Let $A_{c}$ denote the cross-sectional area of the cable and $E_{c}$ denote the modulus of elasticity of the material, then the additional tension in the cable caused by the live load $p$ is
$$
h=\frac{E_{c}A_{c}}{L_{c}}\triangle L_{c}=\frac{E_{c}A_{c}}{L_{c}}\frac{q}{H}\int_0^L w(x)dx.\eqno (8.10)
$$
Substituting (8.10) in (8.5), one obtains
$$
EI w''''(x)-(H+\frac{E_{c}A_{c}}{L_{c}}\frac{q}{H}\int_0^L w(x)dx)w''(x)+\frac{q^{2}}{H^{2}}\frac{E_{c}A_{c}}{L_{c}}\int_0^L w(x)dx=p(x),\ \ \ \ x\in (0,L),\eqno (8.11)
$$
which is the classical form of Melan equation (1.4).

\vskip 3mm


\noindent{\bf Declaration of Interest Statement}

\noindent The author declares that there are no competing interests, financial or otherwise, that could
influence or bias the work reported in this manuscript. The author confirms sole responsibility
for the conception, analysis, and interpretation of the results presented in this study.

\noindent{\bf Acknowledgements}

\noindent This work was supported by the NSFC(No.12261057,\ No.12461035) and the Natural Science Foundation of Gansu Province, China
(grant no.: 22JR5RA264).

\vskip 12mm

\centerline {\bf REFERENCES}\vskip5mm\baselineskip 0.45cm
\begin{description}
\baselineskip 15pt

\item{[1]}~ S. G. Bounopane, D. P. Billington, Theory and history of suspension bridge design from 1823 to 1940, Jour. Strut. Engr. 119 (3) (1993), 954-977.

\item{[2]}~ C. L. Navier, M\'{e}moire sur les ponts suspendus. Imprimerie Royale, Paris, 1823.

\item{[3]}~ W. J. M. Rankine, A manual of applied mechanics. Charles Griffin and Company, London, 1858.

\item{[4]}~ J. Melan, Theory of arches and suspension bridges. Myron Clark Publ. Comp., London, 1913 (German original third edition: Handbuch der Ingenieurwis-senschaften 2, 1906).

\item{[5]}~ D. B. Steinman, A practical treatise on suspension bridges. 2nd Ed., John Wiley, New York, 1929. (Reprinted, 1953.)

\item{[6]}~ D. M. Brotton, A general computer program for the solution of suspension bridge problems, Structural engineering. 44 (1966), 161-167.

\item{[7]}~ S. A. Saafan, Theoretical analysis of suspension bridges, Journal of the structural division. 92 (1966), 1-12.

\item{[8]}~ A. Pugsley, The theory of suspension bridge. Edward Arnold Pub. Ltd, London, 1968. 


\item{[9]}~ M. Irvine, Cable structures. MIT University Press, Cambridge, Mass., 1981. (Reprinted, 1992.)

\item{[10]}~ C. C. Ulstrup. Rating and preliminary analysis of suspension
bridges, Journal of Bridge Engineering, ASCE. 119 (1993), 2653-2679.

\item{[11]}~ G. H. Li, Stability and vibration of bridge. 2nd Ed., China Railway Publishing House, Beijing, 1996. (in Chinese)

\item{[12]}~ P. Walter, Cable-suspended bridges. In: Structural steel designer's handbook. New York: McGraw-Hill Inc, 1999, 15.60-15.67.

\item{[13]}~ P. Clemente, G. Nicolosi, A. Raithel, Preliminary design of very long-span suspension bridges, Engineering Structures. 22 (2000), 1699-1706.

\item{[14]}~ G. P. Wollmann, Preliminary analysis of suspension bridges, J. Bridge Eng. 6 (2001), 227-233.

\item{[15]}~ D. Cobo del Arco, A. C. Aparicio, Preliminary static analysis of suspension bridges, Eng. Struct. 23 (2001), 1096-1103.

\item{[16]}~ J. Q. Lei, M. Z. Zheng, G. Y. Xu, Design of suspension bridge. China Communications Press, Beijing, 2002. (in Chinese)

\item{[17]}~ Y. G. Tan, Z. Zhang, Y. Xu, X. Zhao, Preliminary static analysis of self-anchored suspension bridges, Structural engineering and mechanics, 34 (2010), 281-284.

\item{[18]}~ D. H. Choi, S. G. Gwon, H. Yoo, Nonlinear static analysis of continuous multi-span suspension bridges, International Journal of Steel Structures. 13 (2013), 103-115.

\item{[19]}~ S. U. Shin, M. R. Jung, J. Park, M. Y. Kim, A deflection theory and its validation of earth-anchored suspension bridges under live loads, Ksce Journal of Civil Engineering, 19(1) (2015), 200-212.

\item{[20]}~ M. R. Jung, S. U. Shin, M. M. Attard, M. Y. Kim, Deflection theory for self-anchored suspension bridges under Live Load, Journal of Bridge Engineering, 20 (2015): 04014093.

\item{[21]}~ T. V. Karman, M. A. Biot, Mathematical methods in engineering: an introduction to the mathematical treatment of engineering problems. McGraw-Hill, New York, 1940.

\item{[22]}~ S. P. Timoshenko, D. H. Young, Theory of structures. McGraw Hill, New York, 1968.

\item{[23]}~ M. T. Godard, Recherches sur le Calcul de la resistance des tabliers des Ponts suspendus. Ann. Ponst et Chaussees, 8 (1894), 105-189.

\item{[24]}~ D. J. Peery, An influence-line analysis for suspension bridge. Transactions of ASCE, 121 (1956), 463-480.

\item{[25]}~ H. Bleich, Die Berechnung Verankerter Hangebr$\ddot{u}$cken,  Berlin, 1935.

\item{[26]}~ A. Jennings, Gravity stiffness of classical suspension bridges, Journal of Structural Engineering. 109 (1983), 16-36.

\item{[27]}~ H. Ohshima, K. Sato, N. Watanabe, Structural analysis of suspension bridges, Jour. Engr. Mech. 110 (3) (1984), 392-404.












\item{[28]}~ B. Semper, Finite element methods for suspension bridge models, Computers. Math. Applic. 25(5) (1993),  77-91.

\item{[29]}~ B. Semper, Finite element approximation of a fourth order integro-differential equation, Appl. Math. Lett. (1) (1994), 59-62.

\item{[30]}~ R. C. Mittal, R. Jiwari, A spectral method for suspension bridge model, International Journal of Applied Mathematics and Mechanics. 5(5) (2009), 66-75.

\item{[31]}~ Q. W. Ren, Q. Q. Zhuang, Legendre-galerkin spectral approximation of a class of fourth-order integro-differential equations, Mathematica Numerica Sinica. 35(2) (2013), 125-136.

\item{[32]}~ E. Aruchunan, Y. Wu, B. Wiwatanapataphee, et al. A new variant of arithmetic mean iterative method for fourth order integro-differential equations solution, 2015 IEEE Third International Conference on Artificial Intelligence, Modelling and Simulation, IEEE. 2015.

\item{[33]}~ A. Zeeshan, M. Atlas, Optimal solution of integro-differential equation of suspension bridge model using Genetic Algorithm and Nelder-Mead method, Journal of the Association of Arab Universities for Basic and Applied Sciences. 24(1) (2017), 310-314.

\item{[34]}~ A. C. Lazer, P. J. Mckenna, Large-amplitude periodic oscillations in suspension bridge: some new connections with nonlinear analysis, SIAM Rev. 32 (1990), 537-578.

\item{[35]}~ R. H. Plaut, F. M. Davis, Sudden lateral asymmetry and torsional oscillations of section models of suspension bridges, J. Sound and Vibration. 307 (2007), 894-905.

\item{[36]}~ F. Gazzola, R. Pavani, Wide oscillations finite time blow up for solutions to nonlinear fourth order differential equations, Arch. Ration. Mech. Anal. 207 (2013), 717-752.

\item{[37]}~ W. Lacarbonara, Nonlinear Structural Mechanics: Theory, Dynamical Phenomena and Modeling. Springer Publishing Company, Incorporated, 2013.  Springer, New York, 2013.

\item{[38]}~ F. Gazzola, M. Jleli, B. Samet, On the melan equation for suspension bridges, J Fixed Point Theory Appl. 16 (2014), 159188.

\item{[39]}~ H. Amann, Fixed point equations and nonlinear eigenvalue problems in ordered Banach spaces,  SIAM Rev. 18 (1976), 620-709.
\item{[40]}~ G. S. Ladde, V. Lakshmikantham, A. S. Vatsala, Monotone iterative techniques for nonlinear differential equations. Pitman, Boston, 1985.
\item{[41]}~ C. V. Pao, Nonlinear parabolic and elliptic equations. Plenum Press, New York, 1992.
\item{[42]}~ J. J. Nieto, An abstract monotone iterative technique, Nonlinear Anal. 28 (1997), 1923-1933.
\item{[43]}~ I. Rachunkova, Upper and lower solutions and
topological degree, J. Math. Anal. Appl.  234 (1999), 311-327.
\item{[44]}~ C. D. Coster, P. Habets, The lower and upper solutions method for boundary value problems, in: A. Canada, P. Drabek and A. Fonda (Eds.), Handbook of Differential Equations-Ordinary Differential Equations, 2004.
\item{[45]}~ A. Cabada, J. A. Cid, L. Sanchez, Positivity and lower and upper solutions for fourth order boundary value problems, Nonlinear Anal. 67 (2007), 1599-1612.
\item{[46]}~ R. Y. Ma, J. H. Zhang, S. M. Fu, The method of lower and upper solutions for fourth-order two-point boundary value problems, J. Math. Anal. Appl. 215 (1997), 415-422.

\item{[47]}~ R. Y. Ma, J. X. Wang, Y. Long, Lower and upper solution method for the problem of elastic beam, J. Fixed Point Theory Appl. 20 (2018), 1-13.
\item{[48]}~ Y. X. Li, A monotone iterative technique for solving the bending elastic beam equations, Appl. Math. Comput. 217 (2010), 2200-2208.

\item{[49]}~ D. X. Ma, X. Z. Yang, Upper and lower solution method for fourth-order four-point
boundary value problems, J. Comput. Appl. Math. 223 (2009), 543-551.

\item{[50]}~ Y. M. Wang, Monotone iterative technique for numerical solutions of fourth-order
nonlinear elliptic boundary value problems, Appl. Numer. Math. 57 (2007), 1081-1096.

\item{[51]}~ Z. B. Bai, The method of lower and upper solutions for a bending of an elastic beam equation, J. Math. Anal. Appl. 248
(2000), 195-202.


\item{[52]}~ D. Q. Jiang, W. J. Gao, A. Y. Wan, A monotone method for constructing extremal solutions to fourth-order periodic
boundary value problems, Appl. Math. Comput. 132 (2002), 411-421.
\item{[53]}~ Q. Zhang, S. H. Chen, J. H. Lu,  Upper and lower solution method for fourth-order four-point boundary value problems, J. Comput. Appl. Math. 196 (2006), 387-393.

\item{[54]}~ J. Ehme, P. W. Eloe, J. Henderson,   Upper and lower solution methods for fully nonlinear
boundary value problems, J. Differ. Equations. 180 (2002), 51-64.

\item{[55]}~ P. Habets, L. Sanchez, A monotone method for fourth order boundary value problems involving a factorizable linear
operator, Port. Math. 64(3) (2007), 255-279.

\item{[56]}~ F. Minhos, T. Gyulov, A. I. Santos, Existence and location result for a fourth order boundary value
problem, Discrete Contin. Dyn. Syst. Suppl. (2005), 662-671.

\item{[57]}~ R. Vrabel, On the lower and upper solutions method for the problem of elastic beam with hinged ends. J. Math. Anal. Appl. 421(2) (2015), 1455-1468.

\item{[58]}~ E. Alves, T. F. Ma, L. P. Mauricio, Monotone positive solutions for a fourth order equation with nonlinear boundary conditions, Elsevier Science Publishers B. V. 2009.

\item{[59]}~ V. Berinde,  Iterative Approximation of Fixed Points. Springer, 2007.

\item{[60]}~ M. Asadi, H. Soleimani, S. M. Vaezpour, B. E. Rhoades, On T-stability of Picard iteration in cone metric spaces, Fixed Point Theory and Applications. 1 (2009), 751090.

\item{[61]}~ Arie Romeijn, Reza Sarkhosh, D. Van Goolen, Parametric study on static behaviour of Self-anchored suspension bridges, Steel Structures. 8 (2008), 91-108.

\end{description}
\end{document}